\documentclass[11pt]{article}
\usepackage{mathrsfs}
\usepackage{bbm}
\usepackage{amsfonts}
\usepackage{amsmath}
\usepackage{latexsym}
\newcommand\bgeq{\begin{equation}}
\newcommand\edeq{\end{equation}}
\newcommand\bgar{\begin{array}}
\newcommand\edar{\end{array}}

\title{Linearized stability for a multi-dimensional free boundary problem
  modeling two-phase tumor growth\footnote{This work is supported by the
  China National Natural Science Fundation under grant number 11171357.}}
\author{Shangbin Cui\footnote{E-mail:  cuishb@mail.sysu.edu.cn}}
\date{{\small Department of Mathematics, Sun Yat-Sen University,
   Guangzhou, Guangdong 510275,} \\ [-0.05cm]
  {\small People's Republic of China}}

\begin{document}

\maketitle

\begin{abstract}
  This paper is concerned with a multi-dimensional free boundary problem modeling
  the growth of a tumor with two species of cells: proliferating cells and quiescent
  cells. This free boundary problem has a unique radial stationary solution. By
  using the Fourier expansion of functions on unit sphere via spherical harmonics,
  we establish some decay estimates for the solution of the linearized system of this
  tumor model at the radial stationary solution, so that proving that the radial
  stationary solution is linearly asymptotically stable when neglecting translations.


\end{abstract}
\section{Introduction}
\setcounter{equation}{0}

\hskip 2em
  Since Greenspan first used free boundary problems of partial differential equations
  to model the growth of solid tumors in 1972 (cf. \cite{Green1, Green2}), many different
  tumor models in terms of free boundary problems of partial differential equations have
  been proposed by different groups of researchers, cf., the reviewing articles
  \cite{tumrev1, tumrev2, Fried1, Fried2, Fried3, tumrev3} and references cited therein.
  Rigorous mathematical analysis of such models has made great progress during the
  past twenty more years, and many interesting results have been obtained, cf.,
  \cite{ChenCuiF, Cui1, Cui2, Cui3, Cui4, CuiEsc1, CuiEsc2, CuiEsc3, CuiFri1, CuiWei, FriHu1,
  FriHu2, FriHu3, FriRei1, FriRei2, WuCui, WZ} and references cited therein. From these
  mentioned references it can be seen that dynamics of such tumor models are usually very
  rich, and some recently developed mathematical tools have played important role in their
  rigorous analysis. Major concern in this topic is the asymptotic behavior of solutions
  as time goes to infinity.

  Models describing the growth of tumors possessing homogeneous structures or consisting
  of only one species of cells usually can be reduced into differential equations in
  Banach spaces possessing parabolic structures. It turns out that asymptotic behavior
  of their solutions can be well treated by using the abstract theory for parabolic
  differential equations in Banach spaces, cf. \cite{Cui3, CuiEsc2, CuiEsc3, WuCui, WZ} and
  references therein (see also \cite{FriHu1, FriHu2, FriHu3, FriRei1, FriRei2} and
  references therein for the analysis by using some classical methods). Unlike this, for
  tumors with inhomogeneous structures, their growth models are more complex and the
  corresponding rigorous analysis is much more difficult. Here we particularly mention
  the models reviewed in \cite{Fried1, Fried2, Fried3}, which describes the growth of
  tumors consisting of more than one species of cells. For those models, some interesting
  results have been obtained in the spherically symmetric case (cf.  \cite{ChenCuiF,
  Cui1, Cui2, Cui4, CuiFri1, CuiWei}). For the general non-symmetric case, however, only
  local well-posedness is established (cf. \cite{ChenFri}), and large-time behavior of
  the solution is totally unclear up to the present.

  In this paper we study a tumor model describing the growth of spherically
  non-symmetric tumors consisting of two species of cells: proliferating cells
  and quiescent cells. Mathematical formulation of this model is the following
  multi-dimensional free boundary problem:
\begin{equation}
  \Delta c=F(c) \quad \mbox{for}\;\; x\in\Omega(t),\;\; t>0,
\end{equation}
\begin{equation}
  c=1 \quad \mbox{for}\;\; x\in\partial\Omega(t), \;\; t>0,
\end{equation}
\begin{equation}
  {\partial p\over\partial t}+\nabla\cdot(\vec{v}p)=
  [K_B(c)-K_Q(c)]p+K_P(c)q \quad \mbox{for}\;\; x\in\Omega(t),\;\; t>0,
\end{equation}
\begin{equation}
  {\partial q\over\partial t}+\nabla\cdot(\vec{v}q)=
  K_Q(c)p-[K_P(c)+K_D(c)]q \quad \mbox{for}\;\; x\in\Omega(t),\;\; t>0,
\end{equation}
\begin{equation}
  p+q=1 \quad \mbox{for}\;\; x\in\Omega(t),\;\; t>0,
\end{equation}
\begin{equation}
  \vec{v}=-\nabla\varpi \quad \mbox{for}\;\; x\in\Omega(t),\;\; t>0,
\end{equation}
\begin{equation}
   \varpi= \gamma\kappa \quad \mbox{for}\;\; x\in\partial\Omega(t),\;\; t>0,
\end{equation}
\begin{equation}
  V_n=\vec{v}\cdot\vec{n} \quad \mbox{for}\;\; x\in\partial\Omega(t),\;\; t>0.
\end{equation}
  Here $\Omega(t)$ is the domain occupied by the tumor at time $t$, $c=c(x,t)$,
  $p=p(x,t)$ and $q=q(x,t)$ are the concentration of nutrient, the density of
  proliferating cells and the density of quiescent cells, respectively, $\vec{v}
  =\vec{v}(x,t)$ is the velocity of tumor cell movement, $\varpi=\varpi(x,t)$
  is the pressure distribution in the tumor, $\kappa$ is the mean curvature of
  the tumor surface whose sign is designated by the convention that $\kappa\geq
  0$ at points where $\partial\Omega(t)$ is convex, $\vec{n}$ is the unit outward
  normal vector of $\partial\Omega(t)$, and $V_n$ is the normal velocity of the
  tumor surface. Besides, $F(c)$ is the consumption
  rate of nutrient by tumor cells, $K_B(c)$ is the birth rate of tumor cells,
  $K_P(c)$ and $K_Q(c)$ are transferring rate of tumor cells from quiescent state
  to proliferating state and from proliferating state to quiescent state,
  respectively, and $K_D(c)$ is the death rate of quiescent cells.
  Typically we have (cf. \cite{PPM})
\begin{equation}
  F(c)=\lambda c,
\end{equation}
\begin{equation}
  K_B(c)=k_Bc, \quad K_D(c)=k_D(1-c), \quad K_P(c)=k_Pc, \quad K_Q(c)=k_Q(1-c),
\end{equation}
  where $\lambda$, $k_B$, $k_D$, $k_P$ and $k_Q$ are positive constants. Finally,
  $\gamma$ is a positive constant and is referred as surface tension coefficient.
  For illustration of biological implications of each equation in the above model,
  we refer the reader to see \cite{Cui4, Fried1, Fried2, Fried3, PPM, TP} and
  references therein.

  By summing up (1.3), (1.4) and using (1.5), we get
\begin{equation}
  \nabla\cdot\vec{v}=K_B(c)p-K_D(c)q \quad \mbox{for}\;\; x\in\Omega(t),\;\; t>0.
\end{equation}
  Besides, we note that due to (1.5), the unknown variables $p$ and $q$ are not
  independent. In what follows we shall keep $p$ only. It follows that the system
  (1.1)--(1.8) reduces into the following one:
\begin{equation}
  \Delta c=F(c) \quad \mbox{for}\;\; x\in\Omega(t),\;\; t>0,
\end{equation}
\begin{equation}
  c=1 \quad \mbox{for}\;\; x\in\partial\Omega(t), \;\; t>0,
\end{equation}
\begin{equation}
  {\partial p\over\partial t}+\vec{v}\cdot\nabla p=f(c,p)
   \quad \mbox{for}\;\; x\in\Omega(t),\;\; t>0,
\end{equation}
\begin{equation}
  \nabla\cdot\vec{v}=g(c,p) \quad \mbox{for}\;\; x\in\Omega(t),\;\; t>0,
\end{equation}
\begin{equation}
  \vec{v}=-\nabla\varpi \quad \mbox{for}\;\; x\in\Omega(t),\;\; t>0,
\end{equation}
\begin{equation}
   \varpi=\gamma\kappa \quad \mbox{for}\;\; x\in\partial\Omega(t),\;\; t>0,
\end{equation}
\begin{equation}
  V_n=-\partial_n\varpi \quad \mbox{for}\;\; x\in\partial\Omega(t),\;\; t>0,
\end{equation}
  where
\begin{eqnarray}
  f(c,p)&=&K_P(c)\!+\big[K_M(c)\!-\!K_N(c)\big]p-\!K_M(c)p^2,\\
  g(c,p)&=&K_M(c)p-K_D(c),
\end{eqnarray}
  and
\begin{equation}
  K_M(c)=K_B(c)+K_D(c), \quad K_N(c)=K_P(c)+K_Q(c).
\end{equation}

  We note that the equations (1.15)--(1.17) can be regarded as an elliptic boundary
  value problem for the unknown $\varpi$. Thus, since $\kappa$ is a (quasi-linear)
  second-order elliptic operator for the unknown function $\rho$ describing the
  free boundary $\partial\Omega(t)$, by solving (1.15)--(1.17) to get $\varpi$ as
  a functional of $c,p,\rho$ and next substituting it into (1.18), we see that the
  equation (1.18) can be reduced into a (quasi-linear) third-order parabolic
  pseudo-differential equation for $\rho$ (containing other unknown functions), cf.
  \cite{Esc, EscSim}. On the other hand, the equation (1.14) is clearly a
  quasi-linear hyperbolic equation for the unknown function $p$ (containing other
  unknown functions). This determines that the above model can neither be treated
  as purely parabolic type equations as in \cite{CuiEsc2, CuiEsc3, WuCui, WZ}, nor
  can it be dealt with as a purely hyperbolic equation, which is the main point
  where the difficulty of the above problem lies.

  As we have mentioned before, local existence and uniqueness of a classical solution
  of the initial value problem for the above system has been well established by
  Chen and Friedman in \cite{ChenFri} in more general setting. In \cite{CuiFri1} and
  \cite{ChenCuiF} it was proved that the above system has a unique radial stationary
  solution under the following general conditions on the given functions $F$, $K_B$,
  $K_D$, $K_P$, and $K_Q$:
\begin{equation}
  F,\;\;K_B,\;\;K_D,\;\;K_P\;\;{\rm and}\;\;K_Q\;\;\mbox{are
  $C^{\infty}$-functions};
\end{equation}
\begin{equation}
  F(0)=0\;\;{\rm and}\;\;  F'(c)>0 \quad {\rm for}\;\;0\leq c\leq 1;
\end{equation}
\begin{equation}
\left\{
\begin{array}{l}
   K_B'(c)>0\;\;\mbox{and}\;\;K_D'(c)<0\;\;{\rm for}\;\;0\leq c\leq 1,\;\;
   K_B(0)=0\;\;{\rm and}\;\;K_D(1)=0;\\
   K_P\;\;\mbox{and}\;\;K_Q\;\;\mbox{satisfy the same conditions as}
   \;\; K_B\;\;\mbox{and}\;\;K_D,\;\;\mbox{respectively};\\
   K_B'(c)+K_D'(c)>0\;\;{\rm for}\;\; 0\leq c\leq 1.
\end{array}
\right.
\end{equation}
  Moreover, in \cite{Cui2} it was proved that this unique radial stationary
  solution is asymptotically stable under radial perturbations (see \cite{Cui4}
  for an extension of this result). Naturally, we want to know if this unique
  radial stationary solution is also asymptotically stable under non-radial
  perturbations? Up to now we are unable to give a satisfactory answer to such a
  difficult question. The purpose of this paper is to prove the following weaker
  result:
\medskip

  {\bf Theorem 1.1}\ \ {\em Let $K_B$, $K_D$, $K_P$ and $K_Q$ be given by
  $(1.10)$ with coefficients satisfying the following conditions:
\begin{equation}
   k_B>k_D\geq 2k_Q>0, \qquad  k_B>k_P, \qquad  k_Bk_Q\leq k_Dk_P.
\end{equation}
  There exists a constant $\gamma^*>0$ such that for $\gamma>\gamma^*$, the
  unique radial stationary solution of the system $(1.12)$--$(1.18)$ is
  linearly asymptotically stable modula translations, i.e., the trivial solution
  of the linearized system of $(1.12)$--$(1.18)$ at its unique radial stationary
  solution is asymptotically stable modula translations in suitable function
  spaces.}
\medskip

  Here the phrase ``modula translations'' is used to refer to the following
  property of the system $(1.12)$--$(1.18)$: Since this system is invariant under
  translations in the coordinate space, its stationary solutions are not isolated,
  and by translating a given stationary solution we get a $n$-parameter family of
  stationary solutions. Thus stationary solutions of the above system form a
  $n$-dimensional manifold --- the so-called {\em center manifold}. It follows
  that the trivial solution of the linearized system is also not an isolated
  stationary solution, but instead, all stationary solutions of the linearized
  system make up a $n$-dimensional linear space. Hence, to study asymptotic stability
  of the trivial solution for the linearized system, we must make analysis in certain
  quotient spaces. See Theorem 8.1 in Section 8 for more explicit statement of
  the above result.

  We remark that in (1.25), the condition $k_B>k_D$ is essential and it cannot
  be removed. Indeed, if this condition is removed then the system (1.12)--(1.18)
  does not have a stationary solution and the tumor will finally disappear, cf.
  \cite{CuiFri2}. Unlike this, the other conditions in (1.25) are imposed just
  for technical reasons; see Lemma 3.2 in Section 3 and Lemmas 6.2--6.4 in
  Section 6. We conjecture that the rest conditions can be removed without
  affecting the validity of the above result.

  Throughout this paper we shall make discussion for the general $n$-dimension
  version of the system (1.12)--(1.18) with $n\geq 2$. This will enable us
  to use some abstract theory of differential equations and spherical harmonic
  functions and avoid using concrete expressions of Bessel functions and
  $3$-dimensional spherical harmonics. We note that all discussions made in the
  literature \cite{ChenFri}, \cite{ChenCuiF}, \cite{Cui2} and \cite{CuiFri1}
  can be easily extended to the general $n$-dimension case, so that, in
  particular, the conditions (1.22)--(1.24) ensure that the system
  $(1.12)$--$(1.18)$ has also a unique radial stationary solution in the
  general $n$-dimension case.

  The structure of the rest part is as follows: In the next section we compute
  the linearization of the system (1.10)--(1.16) around the radial stationary
  solution and use the spherical harmonic expansion to make reduction to the
  linearized system. In section 3 we collect a few preliminary lemmas. In
  Sections 4--7 we step by step establish decay estimates for each mono-mode
  system obtained from the spherical harmonic expansion of the linearized system.
  In the last section we compose up all the mono-mode estimates to get desired
  result.

\section{Linearization and reduction}
\setcounter{equation}{0}

\hskip 2em
  We denote by $(c_s(r),p_s(r),v_s(r),\sigma_s(r),\Omega_s)$ ($\Omega_s=\{x\in\mathbb{R}^n:
  r=|x|<R_s$) the unique radial stationary solution of the system (1.10)--(1.16),
  namely, the solution of the following system of equations:
\begin{equation}
  c_s''(r)+\frac{n\!-\!1}{r}c_s'(r)=F(c_s(r)), \quad 0<r<R_s,
\end{equation}
\begin{equation}
  c_s'(0)=0,\quad c_s(R_s)=1,
\end{equation}
\begin{equation}
  v_s(r)p_s'(r)=f(c_s(r),p_s(r)), \quad 0<r<R_s,
\end{equation}
\begin{equation}
  v_s'(r)+\frac{n\!-\!1}{r}v_s(r)=g(c_s(r),p_s(r)), \quad 0<r<R_s,
\end{equation}
\begin{equation}
  v_s(r)=-\sigma_s'(r), \quad 0<r<R_s,
\end{equation}
\begin{equation}
  v_s(0)=0,\quad v_s(R_s)=0.
\end{equation}
  Later on we shall use the following notations:
\begin{equation*}
\begin{array}{c}
  f^*(r)=f(c_s(r),p_s(r)), \qquad   g^*(r)=f(c_s(r),p_s(r)),
\\
  f_c^*(r)=f_c(c_s(r),p_s(r)), \qquad  f_p^*(r)=f_p(c_s(r),p_s(r)),
\\
  g_c^*(r)=g_c(c_s(r),p_s(r)), \qquad  g_p^*(r)=g_p(c_s(r),p_s(r)).
\end{array}
\end{equation*}
  As we mentioned before, existence and uniqueness of the above system has been
  proved in \cite{CuiFri1, ChenCuiF} in the $3$-dimension case. For the general
  $n$-dimension case ($n\geq 2$), the argument is quite similar so that we omit it
  here. Moreover, this solution satisfies the following properties (cf.
  \cite{CuiFri1}):
\begin{equation}
  0<c_s(r)<1\;\; \mbox{for}\;\; 0\leq r<R_s,  \quad c_s'(r)>0\;\; \mbox{for}\;\; 0<r\leq R_s;
\end{equation}
\begin{equation}
  0<p_s(r)<1\;\; \mbox{for}\;\; 0\leq r<R_s,  \quad p_s'(r)>0\;\; \mbox{for}\;\; 0<r\leq R_s;
\end{equation}
\begin{equation}
  p_s(r)>\alpha(c_s(r))\;\; \mbox{for}\;\; 0<r<R_s, \quad
  p_s(0)=\alpha(c_s(0)) \quad \mbox{and} \quad
  p_s(R_s)=\alpha(c_s(R_s)),
\end{equation}
  where
\begin{equation*}
  \alpha(\lambda)=\frac{1}{2K_M(\lambda)}[K_M(\lambda)-K_N(\lambda)+\sqrt{[K_M(\lambda)-K_N(\lambda)]^2
  +4K_M(\lambda)K_P(\lambda)}]
\end{equation*}
  (for $ 0\leq\lambda\leq 1$), and there exist positive constants $c_1$, $c_2$ such that
\begin{equation}
  -c_1r(R_s-r)\leq v_s(r)\leq-c_2r(R_s-r)\;\; \mbox{for}\;\; 0\leq r\leq R_s.
\end{equation}

  It is easy to verify that the linearization of the system (1.10)--(1.16) around
  the stationary solution $(c_s(r),p_s(r),v_s(r),P_s(r),\Omega_s)$ is as follows:
\begin{eqnarray}
  \Delta\sigma&=&F'(c_s(r))\sigma, \quad x\in\Omega_s,\;\; t>0,
\\[0.1cm]
  \sigma|_{x=R_s\omega} &=&-c_s'(R_s)\eta(\omega,t), \quad
  \omega\in\mathbb{S}^{n-1},\;\; t>0,
\\[0.1cm]
  \varphi_t+v_s(r)\varphi_r&=&p_s'(r)\psi_r+f_c^*(r)\sigma+f_p^*(r)\varphi,
  \quad x\in\Omega_s,\;\; t>0,
\\[0.1cm]
  \vec{w}&=&-\nabla\psi, \quad x\in\Omega_s,\;\; t>0,
\\[0.1cm]
  \Delta\psi &=&-g_c^*(r)\sigma-g_p^*(r)\varphi, \quad x\in\Omega_s,\;\; t>0,
\\[0.1cm]
  \psi|_{x=R_s\omega} &=&-{\gamma\over R_s^2}[\eta(\omega,t)
  +{1\over n\!-\!1}\Delta_\omega\eta(\omega,t)], \quad
  \omega\in\mathbb{S}^{n-1},\;\; t>0,
\\[0.1cm]
  \eta_t(\omega,t)&=&-\psi_r(R_s\omega,t)+g(1,1)\eta(\omega,t),
  \quad \omega\in\mathbb{S}^{n-1},\;\; t>0,
\end{eqnarray}
  where $\sigma=\sigma(x,t)$, $\varphi=\varphi(x,t)$, $\vec{w}=\vec{w}(x,t)$,
  $\psi=\psi(x,t)$ and $\eta=\eta(\omega,t)$ ($x\in\Omega_s$, $\omega\in\mathbb{S}^{n-1}$,
  $t\geq0$) are unknown variables, the subscript $r$ denotes derivative in radial direction
  (e.g., $\varphi_r=\frac{\partial\varphi}{\partial r}=\frac{x}{r}\cdot\nabla\varphi$
  etc.), $\Delta_\omega$ denotes the Laplace-Beltrami operator on the unit
  sphere $\mathbb{S}^{n-1}$. To get the above equations, we let
\begin{eqnarray}
\left\{
\begin{array}{l}
  c(x,t)=c_s(r)+\varepsilon\sigma(x,t), \qquad
  p(x,t)=p_s(r)+\varepsilon\varphi(x,t),\quad
\\[0.1cm]
  \vec{v}(x,t)=v_s(r)\omega+\varepsilon\vec{w}(x,t), \quad
  \varpi(x,t)=\sigma_s(r)+\varepsilon\psi(x,t), \qquad
\\[0.1cm]
  \Omega(t)=\{x\in {\mathbb R}^n:\; r<R_s+\varepsilon\eta(\omega,t)\},
\end{array}
\right.
\end{eqnarray}
  where $r=|x|$, $\omega=x/|x|$, and $\varepsilon$ is a small real parameter. Substituting
  these expressions into (1.12)--(1.18) and using some similar arguments as in
  \cite{CuiEsc2}, we obtain (2.11)--(2.17). As an example we only give the deduction
  of the equation (2.13). Substituting the first three relations in (2.18) into the
  equation (1.14), we get
$$
  \varepsilon\varphi_t+[v_s(r)\omega+\varepsilon\vec{w}]\cdot
  [\nabla p_s(r)+\varepsilon\nabla\varphi]
  =f(c_s(r)+\varepsilon\sigma,p_s(r)+\varepsilon\varphi),
$$
  or
$$
\begin{array}{rl}
  &v_s(r)\omega\cdot\nabla p_s(r)+\varepsilon[\varphi_t+v_s(r)\omega\cdot\nabla\varphi
  +\vec{w}\cdot\nabla p_s(r)]+o(\varepsilon)
\\
   =&f(c_s(r),p_s(r))+\varepsilon [f_c^*(r)\sigma+f_p^*(r)\varphi]
   +o(\varepsilon).
\end{array}
$$
  Since $v_s(r)\omega\cdot\nabla p_s(r)=v_s(r)p_s'(r)=f(c_s(r),p_s(r))$, by first
  removing these terms in the above equation, next dividing both sides with
  $\varepsilon$, and finally letting $\varepsilon\to 0$, we get
\begin{equation*}
  \varphi_t+v_s(r)\omega\cdot\nabla\varphi+\vec{w}\cdot\nabla p_s(r)
  =f_c^*(r)\sigma+f_p^*(r)\varphi.
\end{equation*}
  Since $\omega\cdot\nabla\varphi=\varphi_r$ and $\vec{w}\cdot\nabla p_s(r)=-\nabla\psi\cdot
  p_s'(r)\omega=-p_s'(r)\psi_r$ (by (2.14)), we see that (2.13) follows.

  The system (2.11)--(2.17) can be reduced into a 2-system of linear evolution equations
  containing only the unknowns $\varphi$ and $\eta$. To see this we denote by
  $\mathscr{K}$, $\mathscr{K}_0$ and $\mathscr{G}$ respectively the following operators:
  For $\eta\in C^2(\mathbb{S}^{n-1})$, let $u=\mathscr{K}(\eta)\in C^{2*}(\overline{\Omega}_s)$
  and $v=\mathscr{K}_0(\eta)\in C^{2*}(\overline{\Omega}_s)$, where $C^{2*}(\overline{\Omega}_s)$
  represents the second-order Zygmund space in $\overline{\Omega}_s$, be respectively
  solutions of the following elliptic boundary value problems:
$$
\left\{
\begin{array}{l}
  \Delta u=F'(c_s(r))u, \quad x\in\Omega_s,
\\
  u|_{x=R_s\omega}=\eta(\omega), \quad \omega\in\mathbb{S}^{n-1};
\end{array}
\right.
$$
$$
\left\{
\begin{array}{l}
  \Delta v=0, \quad x\in\Omega_s,
\\
  v|_{x=R_s\omega}=\eta(\omega), \quad \omega\in\mathbb{S}^{n-1}.
\end{array}
\right.
$$
  For $h\in C(\overline{\Omega}_s)$, let $w=\mathscr{G}(\eta)\in C^{2*}(\overline{\Omega}_s)$
  be the solution of the following elliptic boundary value problem:
$$
\left\{
\begin{array}{l}
  \Delta w=h, \quad x\in\Omega_s,
\\
  w=0, \quad   x\in\partial\Omega_s.
\end{array}
\right.
$$
  Then from (2.11) and (2.12) we have
$$
   \sigma=-c_s'(R_s)\mathscr{K}(\eta),
$$
  and from (2.15) and (2.16) we have
$$
  \psi=\Phi+\Upsilon+\Psi,
$$
  where
$$
\left\{
\begin{array}{l}
  \Phi=-\mathscr{G}[g_p^*(r)\varphi],
\\
  \Upsilon=-\mathscr{G}[g_c^*(r)\sigma]
  =c_s'(R_s)\mathscr{G}[g_c^*(r)\mathscr{K}(\eta)],
\\
  \Psi=\displaystyle-{\gamma\over R_s^2}\mathscr{K}_0[\eta(\omega,t)
  +{1\over n\!-\!1}\Delta_\omega\eta(\omega,t)].
\end{array}
\right.
$$
  Substituting these expressions into (2.13) and (2.17), we see that the system
  (2.11)--(2.17) reduces into the following $2$-system:
\begin{eqnarray}
\left\{
\begin{array}{l}
  \partial_t\varphi=\mathscr{A}(\varphi,\eta),
\\
  \partial_t\eta=\mathscr{B}(\varphi,\eta),
\end{array}
\right.
\end{eqnarray}
  where
$$
\begin{array}{rcl}
  \mathscr{A}(\varphi,\eta)&=&-v_s(r)\partial_r\varphi+f_p^*(r)\varphi+p_s'(r)\partial_r\Phi
  +p_s'(r)\partial_r\Upsilon+p_s'(r)\partial_r\Psi+f_c^*(r)\sigma
\\
  &=&-v_s(r)\partial_r\varphi+f_p^*(r)\varphi-p_s'(r)\partial_r\mathscr{G}[g_p^*(r)\varphi]
  +c_s'(R_s)p_s'(r)\partial_r\mathscr{G}[g_c^*(r)\mathscr{K}(\eta)]
\\
  &&\displaystyle-{\gamma\over R_s^2}p_s'(r)\partial_r\mathscr{K}_0(\eta
  +{1\over n\!-\!1}\Delta_\omega\eta)
  -c_s'(R_s)f_c^*(r)\mathscr{K}(\eta),
\\
  \mathscr{B}(\varphi,\eta)&=&\displaystyle-\partial_r\Phi|_{r=R_s}
  -\partial_r\Upsilon|_{r=R_s}-\partial_r\Psi|_{r=R_s}+g(1,1)\eta
\\
  &=&\displaystyle\partial_r\mathscr{G}[g_p^*(r)\varphi]|_{r=R_s}
  -c_s'(R_s)\partial_r\mathscr{G}[g_c^*(r)\mathscr{K}(\eta)]|_{r=R_s}
\\
  &&+\displaystyle{\gamma\over R_s^2}\partial_r\mathscr{K}_0(\eta
  +{1\over n\!-\!1}\Delta_\omega\eta)|_{r=R_s}+g(1,1)\eta.
\end{array}
$$

  Let $Y_k(\omega)$, where $\omega$ represents a variable in the sphere $\mathbb{S}^{n-1}$,
  be a spherical harmonics of degree $k$ (cf. \cite{SW}), i.e., $Y_k(\omega)$ is a
  nontrivial solution of the following equation:
$$
  \Delta_\omega Y_k(\omega)=-\lambda_k Y_k(\omega), \quad \mbox{where}\;\;
  \lambda_k=(n+k-2)k
$$
  $(k=0,1,2,\cdots)$. Consider a solution of (2.19) of the form
$$
  \varphi(x,t)=\varphi_k(r,t)Y_k(\omega), \quad
  \eta(\omega,t)=\eta_k(t)Y_k(\omega),
$$
  where $r=|x|$ and $\omega=\displaystyle\frac{x}{|x|}$. Using the identity
$$
  \Delta\varphi=\frac{\partial^2\varphi}{\partial r^2}
  +\frac{n\!-\!1}{r}\frac{\partial\varphi}{\partial r}
  +\frac{1}{r^2}\Delta_{\omega}\varphi
$$
  we easily see that
$$
  \mathscr{K}_0(Y_k)=R_s^{-k}r^kY_k(\omega), \qquad
  \mathscr{K}(Y_k)=R_s^{-k}r^ku_k(r)Y_k(\omega),
$$
  where $u_k$ is the solution of the following problem:
\begin{eqnarray}
\left\{
\begin{array}{l}
  \displaystyle u_k''(r)+\frac{n\!+\!2k\!-\!1}{r}u_k'(r)=F'(c_s(r))u_k(r),
  \quad 0<r<R_s
\\
  u_k'(0)=0,\quad u_k(R_s)=1.
\end{array}
\right.
\end{eqnarray}
  Moreover, for any $f\in C[0,R_s]$ we have
\begin{eqnarray}
  \mathscr{G}(f(r)Y_k(\omega))
  &=&-\Big(r^k\int_r^{R_s}\!\!\!\int_0^{\xi}\frac{\rho^{n+k-1}}{\xi^{n+2k-1}}
  f(\rho)d\rho d\xi\Big)Y_k(\omega)
\nonumber\\
  &=&\frac{1}{n\!+\!2(k\!-\!1)}\Big[\Big(\frac{1}{R_s^{n+2(k-1)}}
  -\frac{1}{r^{n+2(k-1)}}\Big)r^k\int_0^{R_s}\rho^{n+k-1}f(\rho)d\rho
\nonumber\\
  &&+\frac{1}{r^{n+k-2}}\int_r^{R_s}\rho^{n+k-1}f(\rho)d\rho
  -r^k\int_r^{R_s}\rho^{-k+1}f(\rho)d\rho\Big]Y_k(\omega).
\end{eqnarray}
  Hence
\begin{eqnarray*}
  \sigma&=&-c_s'(R_s)R_s^{-k}r^ku_k(r)Y_k(\omega)\eta_k(t),
\\
  \Phi&=&\displaystyle\Big(r^k\int_r^{R_s}\!\!\!\int_0^{\xi}\frac{\rho^{n+k-1}}{\xi^{n+2k-1}}
  g_p^*(\rho)\varphi_k(\rho,t)d\rho d\xi\Big)Y_k(\omega),
\\
  \Upsilon&=&\displaystyle-c_s'(R_s)R_s^{-k}\Big(r^k\int_r^{R_s}\!\!\!\int_0^{\xi}\frac{\rho^{n+2k-1}}{\xi^{n+2k-1}}
  g_c^*(\rho)u_k(\rho)d\rho d\xi\Big)Y_k(\omega)\eta_k(t),
\\
  \Psi&=&\displaystyle-\Big(1-{\lambda_k\over n\!-\!1}\Big)\gamma R_s^{-k-2}r^kY_k(\omega)\eta_k(t).
\end{eqnarray*}
  It follows that
\begin{eqnarray*}
  \mathscr{A}(\varphi,\eta)&=&-v_s(r)\partial_r\varphi+f_p^*(r)\varphi+p_s'(r)\partial_r\Phi
  +p_s'(r)\partial_r\Upsilon+p_s'(r)\partial_r\Psi+f_c^*(r)\sigma
\\
  &=&\displaystyle\Big\{-v_s(r)\partial_r\varphi_k+f_p^*(r)\varphi_k+p_s'(r)\frac{\partial}{\partial r}
  \Big(r^k\int_r^{R_s}\!\!\!\int_0^{\xi}\frac{\rho^{n+k-1}}{\xi^{n+2k-1}}g_p^*(\rho)\varphi_k(\rho,t)d\rho d\xi\Big)
\\
  &&\displaystyle-c_s'(R_s)R_s^{-k}p_s'(r)\frac{\partial}{\partial r}\Big(r^k\int_r^{R_s}\!\!\!\int_0^{\xi}
  \frac{\rho^{n+2k-1}}{\xi^{n+2k-1}}g_c^*(\rho)u_k(\rho)d\rho d\xi\Big)\eta_k(t)
\\
  &&\displaystyle-\Big(1-{\lambda_k\over n\!-\!1}\Big)\gamma kR_s^{-k-2}r^{k-1}p_s'(r)\eta_k(t)
  -c_s'(R_s)R_s^{-k}f_c^*(r)r^ku_k(r)\eta_k(t)\Big\}Y_k(\omega),
\\
  \mathscr{B}(\varphi,\eta)&=&\displaystyle-\partial_r\Phi|_{r=R_s}
  -\partial_r\Upsilon|_{r=R_s}-\partial_r\Psi|_{r=R_s}+g(1,1)\eta
\\
  &=&\displaystyle\Big\{-\frac{\partial}{\partial r}\Big(r^k\int_r^{R_s}\!\!\!\int_0^{\xi}
  \frac{\rho^{n+k-1}}{\xi^{n+2k-1}}g_p^*(\rho)\varphi_k(\rho,t)d\rho d\xi\Big)\Big|_{r=R_s}
\\
  &&\displaystyle+c_s'(R_s)R_s^{-k}\frac{\partial}{\partial r}\Big(r^k\int_r^{R_s}\!\!\!\int_0^{\xi}
  \frac{\rho^{n+2k-1}}{\xi^{n+2k-1}}g_c^*(\rho)u_k(\rho)d\rho d\xi\Big)\Big|_{r=R_s}\eta_k(t)
\\
  &&+\displaystyle\Big(1-{\lambda_k\over n\!-\!1}\Big)\gamma kR_s^{-3}\eta_k(t)
  +g(1,1)\eta_k(t)\Big\}Y_k(\omega).
\end{eqnarray*}
  Hence we get
\begin{eqnarray}
  \frac{\partial\varphi_{k}}{\partial t}&=&\mathscr{L}_k(\varphi_{k})
  +b_k(r,\gamma)\eta_{k},
\\
  \frac{d\eta_{k}}{dt}&=&\alpha_k(\gamma)\eta_{k}+J_k(\varphi_{k}),
\end{eqnarray}
  where
\begin{eqnarray*}
  \alpha_k(\gamma)&=&\Big(1-\frac{\lambda_k}{n\!-\!1}\Big)
  \frac{k\gamma}{R_s^3}+g(1,1)-\frac{c_s'(R_s)}{R_s^{n+2k-1}}
  \int_0^{R_s}\rho^{n+2k-1}g_c^*(\rho)u_k(\rho)d\rho,
\\
  b_k(r,\gamma)&=&\displaystyle-c_s'(R_s)R_s^{-k}p_s'(r)
  \frac{\partial}{\partial r}\Big(r^k\int_r^{R_s}\!\!\!
  \int_0^{\xi}\frac{\rho^{n+2k-1}}{\xi^{n+2k-1}}g_c^*(\rho)u_k(\rho)d\rho d\xi\Big)
  \qquad\qquad\qquad\quad
\\
  &&\displaystyle-\Big(1-{\lambda_k\over n\!-\!1}\Big)\gamma
  kR_s^{-k-2}r^{k-1}p_s'(r)-c_s'(R_s)R_s^{-k}f_c^*(r)r^ku_k(r)
\\
  &=&\displaystyle-\Big(1-{\lambda_k\over n\!-\!1}\Big)\gamma
  kR_s^{-k-2}r^{k-1}p_s'(r)-c_s'(R_s)R_s^{-k}f_c^*(r)r^ku_k(r)
\\
  &&\displaystyle-c_s'(R_s)R_s^{-k}r^{k-1}p_s'(r)\Big[\theta_k
  \int_r^{R_s}\rho g_c^*(\rho)u_k(\rho)d\rho-\frac{1-\theta_k}{r^{n+2(k-1)}}
  \int_0^r\rho^{n+2k-1}g_c^*(\rho)u_k(\rho)d\rho
\\
  &&\displaystyle-\frac{\theta_k}{R_s^{n+2(k-1)}}
  \int_0^{R_s}\rho^{n+2k-1}g_c^*(\rho)u_k(\rho)d\rho\Big],
\end{eqnarray*}
  where $\theta_k=\displaystyle\frac{k}{n\!+\!2(k\!-\!1)}$, and for $\phi=\phi(r)$,
\begin{eqnarray*}
  \mathscr{L}_k(\phi)&=&\displaystyle-v_s(r)\phi'(r)+f_p^*(r)\phi(r)
  +p_s'(r)\frac{\partial}{\partial r}
  \Big(r^k\int_r^{R_s}\!\!\!\int_0^{\xi}\frac{\rho^{n+k-1}}{\xi^{n+2k-1}}
  g_p^*(\rho)\phi(\rho)d\rho d\xi\Big)
\\
  &=&\displaystyle-v_s(r)\phi'(r)+f_p^*(r)\phi(r)+r^{k-1}p_s'(r)
  \Big[{\theta_k}\int_r^{R_s}\rho^{-k+1}g_p^*(\rho)\phi(\rho)d\rho
\\
  &&\displaystyle-\frac{1-\theta_k}{r^{n+2(k-1)}}
  \int_0^r\rho^{n+k-1}g_p^*(\rho)\phi(\rho)d\rho-\frac{\theta_k}{R_s^{n+2(k-1)}}
  \int_0^{R_s}\rho^{n+k-1}g_p^*(\rho)\phi(\rho)d\rho\Big]
\end{eqnarray*}
  and
\begin{eqnarray*}
  J_k(\phi)&=&\frac{1}{R_s^{n+k-1}}\int_0^{R_s}\rho^{n+k-1}g_p^*(\rho)\phi(\rho)d\rho.
\end{eqnarray*}
  Multiplying (2.23) with $R_s^{-(k-1)}r^{k-1}p_s'(r)$ and adding it into (2.22), we
  see that the system (2.22)--(2.23) reduces into the following equivalent one:
\begin{equation}
  \frac{\partial\widetilde{\varphi}_{k}}{\partial t}
  =\tilde{\mathscr{L}}_k(\widetilde{\varphi}_{k})
  +c_k(r)\eta_{k},
\end{equation}
\begin{equation}
  \frac{d\eta_{k}}{dt}=\widetilde{\alpha}_k(\gamma)\eta_{k}
  +J_k(\widetilde{\varphi}_{k}),
\end{equation}
  where $\widetilde{\varphi}_{k}=\varphi_{k}+R_s^{-(k-1)}r^{k-1}p_s'(r)\eta_k(t)$,
\begin{eqnarray}
  \tilde{\mathscr{L}}_k(\phi)&=&\displaystyle \mathscr{L}_k(\phi)
  +R_s^{-(k-1)}r^{k-1}p_s'(r)J_k(\phi)
\nonumber\\ [0.3cm]
  &=&\displaystyle-v_s(r)\phi'(r)+f_p^*(r)\phi(r)
  +r^{k-1}p_s'(r)\Big[{\theta_k}
  \int_r^{R_s}\rho^{-k+1}g_p^*(\rho)\phi(\rho)d\rho
\nonumber\\ [0.3cm]
  &&\displaystyle+\frac{1-\theta_k}{R_s^{n+2(k-1)}}
  \int_0^{R_s}\rho^{n+k-1}g_p^*(\rho)\phi(\rho)d\rho
  -\frac{1-\theta_k}{r^{n+2(k-1)}}
  \int_0^r\rho^{n+k-1}g_p^*(\rho)\phi(\rho)d\rho\Big],\qquad
\end{eqnarray}
\begin{eqnarray}
  c_k(r)&=&b_k(r,\gamma)+\alpha_k(\gamma)R_s^{-(k-1)}r^{k-1}p_s'(r)
  -R_s^{-(k-1)}\tilde{\mathscr{L}}_k[r^{k-1}p_s'(r)]
\nonumber\\ [0.3cm]
  &=&\displaystyle\frac{r^{k-1}}{R_s^{k-1}}\Big\{[g(1,1)-g^*(r)]p_s'(r)
  +\frac{n\!+\!k\!-\!2}{r}f^*(r)+f_c^*(r)
  \Big[c_s'(r)-c_s'(R_s)R_s^{-1}r u_k(r)\Big]
\nonumber\\ [0.3cm]
  &&\displaystyle  -p_s'(r)\Big[{\theta_k}
  \int_r^{R_s}v_k(\rho)d\rho+\frac{1-\theta_k}{R_s^{n+2(k-1)}}
  \int_0^{R_s}\rho^{n+2(k-1)}v_k(\rho)d\rho
\nonumber\\ [0.3cm]
  &&\displaystyle -\frac{1-\theta_k}{r^{n+2(k-1)}}
  \int_0^r\rho^{n+2(k-1)}v_k(\rho)d\rho\Big]\Big\},
\end{eqnarray}
  where $f^*(r)=f(c_s(r),p_s(r))$, $g^*(r)=g(c_s(r),p_s(r))$,
\begin{equation}
  v_k(r)=g_p^*(r)p_s'(r)+c_s'(R_s)R_s^{-1}g_c^*(r)r u_k(r),
\end{equation}
  and
\begin{eqnarray}
  \widetilde{\alpha}_k(\gamma)&=&\alpha_k(\gamma)-R_s^{-(k-1)}J_k(r^{k-1}p_s'(r))
\nonumber\\
  &=&\displaystyle\Big(1-\frac{\lambda_k}{n\!-\!1}\Big)\frac{k\gamma}{R_s^3}
  +g(1,1)-\frac{1}{R_s^{n+2(k-1)}}\int_0^{R_s}\rho^{n+2(k-1)}v_k(\rho)d\rho.
\end{eqnarray}

  A simple computation shows that $u_1(r)=\displaystyle\frac{R_sc_s'(r)}{rc_s'(R_s)}$
  (see the assertion (4) of Lemma 3.3 below), so that $v_1(r)=g_p^*(r)p_s'(r)
  +g_c^*(r)c_s'(r)=\displaystyle\frac{d}{dr}g^*(r)$. Using these facts, one may easily
  check that
\begin{equation}
  c_1(r)\equiv 0, \qquad \widetilde{\alpha}_1(\gamma)=0.
\end{equation}
  This implies that in the case $k=1$, (2.24)--(2.25) has the following
  stationary solution:
\begin{equation}
  \widetilde{\varphi}_1(r,t)\equiv 0, \qquad \eta_1(t)\equiv const..
\end{equation}
  Or equivalently, in the case $k=1$ the system (2.22)--(2.23) has the following
  stationary solutions:
\begin{equation}
  \varphi_1(r,t)=-cp_s'(r), \qquad \eta_1(t)=c,
\end{equation}
  where $c$ is an arbitrary real constant. This means that the system (2.19) has
  infinite many stationary solutions, and all its stationary solutions form a
  $n$-dimensional linear space (cf. (8.2) in Section 8).

\section{Some preliminary lemmas}
\setcounter{equation}{0}

\hskip 2em
  In this section we collect some preliminary lemmas.
\medskip

  {\bf Lemma 3.1}\ \ {\em The following inequalities hold for all $0\leq r\leq 1$:}
\begin{equation}
  f_p^*(r)<0, \quad  f_c^*(r)>0, \quad  g_p^*(r)>0 \quad \mbox{and}\quad  g_c^*(r)>0.
\end{equation}

  {\em Proof}:\ \ From (1.16) and (1.17) we see that
\begin{eqnarray*}
  f_p^*(r)&=&K_M(c_s(r))-K_N(c_s(r))-2K_M(c_s(r))p_s(r),
\nonumber\\
  f_c^*(r)&=&K_P'(c_s(r))+[K_M'(c_s(r))-K_N'(c_s(r))]p_s(r)-K_M'(c_s(r))p_s^2(r),
\nonumber\\
  g_p^*(r)&=&K_M(c_s(r)),
\nonumber\\
  g_c^*(r)&=&K_M'(c_s(r))p_s(r)-K_D'(c_s(r)).
\end{eqnarray*}
  Since $K_M(c)>0$, $K_M'(c)>0$, $K_P'(c)>0$, $K_D'(c)<0$, $K_N'(c)<0$ for $c>0$ and
  $0<p_s(r)\leq 1$ for $0\leq r\leq 1$, the last three inequalities in (3.1) are
  immediate. Next, since $v_s(r)<0$ for $0<r<1$ and $p_s'(r)>0$ for $0<r\leq 1$, we
  see that $f(c_s(r),p_s(r))=v_s(r)p_s'(r)<0$ for $0<r<1$, which implies that
$$
  p_s(r)\geq\frac{K_M(c_s(r))-K_N(c_s(r))+\sqrt{[K_M(c_s(r))-K_N(c_s(r))]^2
  +4K_M(c_s(r))K_P(c_s(r))}}{2K_M(c_s(r))}
$$
  for $0\leq r\leq 1$. Hence
$$
  f_p^*(r)\leq-\sqrt{[K_M(c_s(r))-K_N(c_s(r))]^2+4K_M(c_s(r))K_P(c_s(r))}<0
$$

  {\bf Lemma 3.2}\ \ {\em Let the conditions in $(1.25)$ be satisfied. There exists
  a constant $c_0>0$ such that as $r\to 0^+$,}
\begin{equation}
  p_s'(r)=\big(c_0+o(1)\big)r.
\end{equation}

  {\em Proof}:\ \ Let $\theta=f_p^*(0)/v_s'(0)$. Since $v_s'(0)<0$ and
  $f_p^*(0)<0$, we have $\theta>0$. By Lemma 5.2 of \cite{ChenCuiF} we know that
  there exists a constant $c_0>0$ such that as $r\to 0^+$,
$$
  p_s'(r)=\big(c_0+o(1)\big)\left\{
\begin{array}{ll}
    r &\quad {\it if}\;\;\kappa>2,\\
    \displaystyle r|\ln r| &\quad {\it if}\;\;\kappa=2,\\
    r^{\theta-1} &\quad {\it if}\;\;\kappa<2.
\end{array}
  \right.
$$
  Hence, we only need to prove that $\theta>2$, or equivalently, $f_p^*(0)<2v_s'(0)$.
  We have
$$
  f_p^*(0)=K_M(c_s(0))-K_N(c_s(0))-2K_M(c_s(0))p_s(0),
$$
$$
  v_s'(0)=\frac{1}{n}g^*(0)=\frac{1}{n}[K_M(c_s(0))p_s(0)-K_D(c_s(0))].
$$
  Hence $\theta>2$ if and only if
\begin{equation}
  n[K_M(c_s(0))-K_N(c_s(0))]+2K_D(c_s(0))<2(n+1)K_M(c_s(0))p_s(0).
\end{equation}
  Since (see (2.9))
$$
  p_s(0)=\frac{K_M(c_s(0))-K_N(c_s(0))+\sqrt{[K_M(c_s(0))-K_N(c_s(0))]^2
  +4K_M(c_s(0))K_P(c_s(0))}}{2K_M(c_s(0))},
$$
  we see that (3.3) is equivalent to
\begin{eqnarray*}
  &&K_M(c_s(0))-K_N(c_s(0))+(n+1)\sqrt{[K_M(c_s(0))-K_N(c_s(0))]^2
  +4K_M(c_s(0))K_P(c_s(0))}
\\
  &>&2K_D(c_s(0)).
\end{eqnarray*}
  This is equivalent to
\begin{eqnarray*}
  &&n(n+2)[K_M(c_s(0))-K_N(c_s(0))]^2+4(n+1)^2K_M(c_s(0))K_P(c_s(0))
\\
  &&+4K_B(c_s(0))K_D(c_s(0))>4K_D(c_s(0))K_N(c_s(0)).
\end{eqnarray*}
  It is easy to check that the conditions in (1.25) ensure that the above inequality
  holds. Hence the desired assertion follows. $\qquad\Box$
\medskip

  {\bf Lemma 3.3}\ \ {\em For the solution $u_k$ of the problem $(2.20)$, we have
  the following assertions:

  $(1)$\ \ $u_k\in C^{\infty}[0,R_s]$, and $0<u_k(r)\leq 1$ for $0\leq r\leq R_s$.

  $(2)$\ \ There exists a constant $C>0$ independent of $k$ such that
\begin{equation}
  1-\frac{C}{n+2k}(R_s-r)\leq u_k(r)\leq 1 \quad \mbox{for}\;\;
  0\leq r\leq R_s,
\end{equation}
\begin{equation}
  0\leq u_k'(r)\leq\frac{Cr}{n+2k} \quad \mbox{for}\;\;
  0\leq r\leq R_s.
\end{equation}

  $(3)$\ \ $u_k(r)$ is monotone non-decreasing in $k$, i.e., $u_k(r)\geq u_l(r)$ for
  $0\leq r\leq R_s$ and $k>l$.

  $(4)$\ \ $u_1(r)=\displaystyle\frac{R_sc_s'(r)}{rc_s'(R_s)}$.

  $(5)$\ \ $u_0(r)>\displaystyle\frac{rc_s'(r)}{R_sc_s'(R_s)}$ for $0\leq r<R_s$.}
\medskip

  {\em Proof}:\ \ The problem (2.20) can be regarded as the spherically symmetric form
  of the $n\!+\!2k$-dimensional elliptic boundary value problem
$$
\left\{
\begin{array}{l}
  \Delta u(x)=F'(c_s(r))u(x) \quad \mbox{for}\;\; |x|<R_s,
\\
  u(x)=1 \quad \mbox{for}\;\; |x|=R_s.
\end{array}
\right.
$$
  From this fact the assertion (1) immediately follow. Next, from the first equation in
  (2.20) we have
$$
  u_k'(r)=\frac{1}{r^{n+2k-1}}\int_0^r\rho^{n+2k-1}F'(c_s(\rho))u_k(\rho)d\rho.
$$
  Let $C_0=\displaystyle\max_{0\leq c\leq 1}F'(c)$. Then we get
$$
  0\leq u_k'(r)\leq\frac{C_0}{r^{n+2k-1}}\int_0^r\rho^{n+2k-1}d\rho
  \leq\frac{C_0r}{n+2k} \quad \mbox{for}\;\; 0\leq r\leq R_s.
$$
  This proves (3.5). Since $u_k(R_s)=1$, by integrating (3.5) over $(r,R_s)$ we get
  (3.4). This proves the assertion (2). The assertion (3) follows from the fact that
  $u_k'(r)\geq 0$ and the maximum principle for second-order elliptic equations.
  The assertion (4) follows from direct computation. Indeed, a direct computation
  shows that the function $\tilde{u}_1(r)=\displaystyle\frac{R_sc_s'(r)}{rc_s'(R_s)}$
  satisfies the same equation as $u_1(r)$ in the region $0<r<R_s$, and it is clear
  that $\tilde{u}_1(R_s)=1$. Since
$$
  \lim_{r\to 0^+}\Big(\frac{c_s'(r)}{r}\Big)'
  =\lim_{r\to 0^+}\frac{rc_s''(r)-c_s'(r)}{r^2}
  =\frac{1}{2}c_s'''(0)=0,
$$
  we see that $\tilde{u}_1'(0)=0$. Hence, by uniqueness of the solution of the elliptic
  boundary value problem we get $\tilde{u}_1(r)=u_1(r)$ for $0\leq r\leq R_s$.
  Finally, it is easy to check that the function $\tilde{u}_0(r)=\displaystyle
  \frac{rc_s'(r)}{R_sc_s'(R_s)}$ satisfies the inequality
$$
  \tilde{u}_0''(r)+\frac{n-3}{r}\tilde{u}_0'(r)
  \geq F'(c_s(r))\tilde{u}_0(r) \quad \mbox{for}\;\; 0<r<R_s.
$$
  Using the fact that $u_0'(r)\geq 0$ we can also easily see that $u_0(r)$ satisfies
  the inequality
$$
  u_0''(r)+\frac{n-3}{r}u_0'(r)
  \leq F'(c_s(r))u_0(r) \quad \mbox{for}\;\; 0<r<R_s.
$$
  Since $u_0'(0)=\tilde{u}_0'(0)$ and $u_0(R_s)=\tilde{u}_0(R_s)$, by the maximum
  principle we see that the desired assertion follows. This completes the proof of
  Lemma 3.3. $\quad\Box$
\medskip

  {\bf Corollary 3.4}\ \ {\em For the function $v_k$ given by $(2.28)$, there exists
  a positive constant $C$ independent of $k$ such that $0\leq v_k(r)\leq C[1+p_s'(r)]$
  for $0\leq r\leq R_s$.}
\medskip

\section{Decay estimates for some positive semigroups}
\setcounter{equation}{0}

\hskip 2em
  In this preliminary section we establish decay estimates for some positive
  semigroups in $C[0,R_s]$ and $L^1([0,R_s],r^{n-1}dr)$. Let $v_s=v_s(r)$ be as before
  and $a=a(r)$ be a real-valued continuous function defined in $[0,R_s]$. Let $L_0$ be
  the following differential operator in $[0,R_s]$: For any function $\varphi$ defined
  in $[0,R_s]$ such that the right-hand side of the following equality makes sense,
$$
  L_0\varphi(r)=-v_s(r)\varphi'(r)+a(r)\varphi(r) \quad
  \mbox{for}\;\; 0\leq r\leq R_s.
$$
  We shall regard $L_0$ both as an unbounded closed linear operator in $C[0,R_s]$ with
  domain $D(L_0)=\{\varphi\in C[0,R_s]\cap C^1(0,R_s):v_s(r)\varphi'(r)\in C[0,R_s]\}$
  and as an unbounded closed linear operator in $L^1([0,R_s],r^{n-1}dr)$ with domain
  $D(L_0)=\{\varphi\in L^1([0,R_s],r^{n-1}dr)\cap W^{1,1}_{\rm loc}(0,R_s):v_s(r)
  \varphi'(r)\in L^1([0,R_s],r^{n-1}dr)\}$. We denote
$$
  \lambda_0(a)=a(R_s), \qquad
  \lambda_1(a)=\max\{a(0),a(R_s)\}, \qquad
  \lambda^*(a)=\displaystyle\max_{0\leq r\leq R_s}a(r).
$$
  It is clear that $\lambda_0(a)\leq\lambda_1(a)\leq\lambda^*(a)$. In what follows the
  notation $\lambda$ denotes a complex number, and $h$ denotes a complex-valued function
  defined in $[0,R_s]$.
\medskip

  {\bf Lemma 4.1}\ \ {\em We have the following assertions:

  $(1)$\ \ If ${\rm Re}\lambda>\lambda_0(a)$ then for any $h\in C(0,R_s]$ the equation
\begin{equation}
  \lambda\varphi(r)-L_0\varphi(r)=h(r) \quad \mbox{for}\;\;0<r<R_s
\end{equation}
  has a unique solution $\varphi=\varphi_{\lambda}\in C(0,R_s]\cap C^1(0,R_s)$, with boundary
  value
\begin{equation}
  \varphi_{\lambda}(R_s)=\frac{h(R_s)}{\lambda-a(R_s)}.
\end{equation}
  Moreover, for any $0<r_0<R_s$ there exists a corresponding constant $C_{\lambda}(r_0)>0$ such that
\begin{equation}
  \max_{r_0\leq r\leq R_s}|\varphi_{\lambda}(r)|\leq C_{\lambda}(r_0)
  \max_{r_0\leq r\leq R_s}|h(r)|.
\end{equation}
  If furthermore $\lambda$ is real and $h(r)\geq 0$ for $0<r\leq R_s$ then also
  $\varphi_{\lambda}(r)\geq 0$ for $0<r\leq R_s$.

  $(2)$\ \ If ${\rm Re}\lambda>\lambda_1(a)$ then for any $h\in C[0,R_s]$ the unique
  solution of $(4.1)$ ensured by the above assertion belongs to $C[0,R_s]\cap
  C^1(0,R_s)$, and in addition to $(4.2)$ we have also that
\begin{equation}
  \varphi_{\lambda}(0)=\frac{h(0)}{\lambda-a(0)}.
\end{equation}
  Moreover, there exists a constant $C_{\lambda}>0$ such that
\begin{equation}
  \max_{0\leq r\leq R_s}|\varphi_{\lambda}(r)|\leq C_{\lambda}
  \max_{0\leq r\leq R_s}|h(r)|.
\end{equation}

  $(3)$\ \ If ${\rm Re}\lambda>\lambda^*(a)$ then the estimate $(4.5)$ can be improved
  as follows:
\begin{equation}
  \max_{0\leq r\leq R_s}|\varphi_{\lambda}(r)|\leq [{\rm Re}\lambda-\lambda^*(a)]^{-1}
  \max_{0\leq r\leq R_s}|h(r)|.
\end{equation}
}

  {\em Proof}:\ \ We first assume that ${\rm Re}\lambda>\lambda_0(a)$. Choose a number
  $r_0\in (0,R_s)$ and set
$$
   W_{\lambda}(r)=\exp\Big(\int^r_{r_0}{\lambda-a(\rho)\over v_s(\rho)}d\rho\Big)
   \quad \mbox{for}\quad 0<r<R_s.
$$
  It is easy to see that $W_{\lambda}\in C^1(0,R_s)$, and
\begin{eqnarray}
  W_{\lambda}(r)&=&C(R_s-r)^{\alpha_1}\big(1+o(1)\big)\;\;\; {\rm as}\;\; r\to R_s^-,
\end{eqnarray}
  where $\alpha_1=\displaystyle{\lambda-a(R_s)\over v_s'(R_s)}$, and $C$ is a nonzero
  constant (depending on the choice of $r_0$). Note that ${\rm Re}\,\alpha_1>0$. Clearly,
  the equation (4.1) can be rewritten as follows:
$$
  \frac{d}{dr}\Big(W_{\lambda}(r)\varphi(r)\Big)={h(r)W_{\lambda}(r)\over v_s(r)}.
$$
 Letting $c=\varphi(r_0)$ and integrating both sides of this equation from $r_0$ to an
 arbitrary point $0<r<R_s$, we see that the general solution of the equation (4.1) is
 given by
\begin{equation}
  \varphi(r)={1\over W_{\lambda}(r)}\Big[c+\int_{r_0}^r
  {h(\eta)W_{\lambda}(\eta)\over v_s(\eta)}d\eta\Big]
  ={1\over W_{\lambda}(r)}\Big[c-\int_{r_0}^r
  {h(\eta)W_{\lambda}(\eta)\over |v_s(\eta)|}d\eta\Big]
\end{equation}
  (for $0<r<R_s$). Since $\displaystyle\lim_{r\to R_s^-}W_{\lambda}(r)=0$ and
$$
  \lim_{r\to R_s^-}\int_{r_0}^r{h(\eta)W_{\lambda}(\eta)\over v_s(\eta)}
  d\eta=\int_{r_0}^{R_s}{h(\eta)W_{\lambda}(\eta)\over v_s(\eta)}d\eta
$$
  is a finite number (by (4.7) and (2.10)), we see that $\displaystyle\lim_{r\to R_s^-}
  \varphi(r)=\infty$ unless $c=\displaystyle\int_{r_0}^{R_s}{h(\eta)W_{\lambda}(\eta)
  \over |v_s(\eta)|}d\eta$, in which case
$$
  \varphi(r)=\frac{h(R_s)}{v_s'(R_s)}\frac{1}{\alpha_1}[1+o(1)]
  =\frac{h(R_s)[1+o(1)]}{\lambda-f_p^*(R_s)} \quad \mbox{as}\;\; r\to R_s^-.
$$
  Hence the solution which is bounded near $r=R_s$ is unique, and this unique bounded
  solution is given by
\begin{equation}
  \varphi_{\lambda}(r)={1\over W_{\lambda}(r)}\int^{R_s}_r{h(\eta)W_{\lambda}(\eta)
  \over |v_s(\eta)|}d\eta  \quad \mbox{for}\;\;0<r<R_s,
\end{equation}
  which is continuous in $(0,R_s]$, continuously differentiable in $(0,R_s)$, and
  satisfies (4.2) (as we have seen above). The estimate (4.3) easily follows
  from (4.9) because the function $r\mapsto\displaystyle{1\over |W_{\lambda}(r)|}
  \int^{R_s}_r\!\!{|W_{\lambda}(\eta)|\over |v_s(\eta)|}d\eta$ is continuous in
  $(0,R_s]$. Moreover, if $\lambda$ is real then $W_{\lambda}(r)>0$ for $0<r<R_s$.
  Using this fact and the expression (4.9) we easily see that if $h\geq 0$ then
  also $\varphi_{\lambda}\geq 0$. This proves the assertion (1).

  Next we assume that ${\rm Re}\lambda>\lambda_1(a)$. Then in addition to (4.7) we
  have also that
\begin{eqnarray}
  W_{\lambda}(r)&=&Cr^{-\alpha_0}\big(1+o(1)\big)\;\;\; {\rm as}\;\; r\to 0^+,
\end{eqnarray}
  where $\alpha_0=\displaystyle\frac{a(0)-\lambda}{v_s'(0)}$,
  and $C$ is a nonzero constant (depending on the choice of $r_0$). Note that
  ${\rm Re}\,\alpha_0>0$. Using this fact we can easily deduce that for any
  constant $c$ the function $\varphi$ given by (4.8) satisfies
$$
  \varphi(r)=-\frac{h(0)}{v_s'(0)}\frac{1}{\alpha_0}[1+o(1)]
  =\frac{h(0)[1+o(1)]}{\lambda-f_p^*(0)} \quad \mbox{as}\;\; r\to 0^+.
$$
  Hence, all solutions given by (4.8) are continuous at $r=0$ and satisfy
  (4.4). In particular, the solution $\varphi_{\lambda}$ given by (4.9)
  belongs to $C[0,R_s]\cap C^1(0,R_s)$ and satisfy both (4.2) and (4.4).
  The estimate (4.5) easily follows from (4.9) because in the present
  situation the function $r\mapsto\displaystyle{1\over |W_{\lambda}(r)|}
  \int^{R_s}_r\!\!{|W_{\lambda}(\eta)|\over |v_s(\eta)|}d\eta$ is continuous
  in $[0,R_s]$. This proves the assertion (2).

  Finally we assume that ${\rm Re}\lambda>\lambda^*(a)$. Recalling $\lambda^*(a)=
  \displaystyle\max_{0\leq r\leq R_s}a(r)$, we have
\begin{eqnarray*}
  \displaystyle{1\over |W_{\lambda}(r)|}\int^{R_s}_r{|W_{\lambda}(\eta)|
  \over |v_s(\eta)|}d\eta
  &\leq&\displaystyle\int^{R_s}_r|v_s(\eta)|^{-1}e^{-\int^{\eta}_r{{\rm Re}
  \lambda-a(\rho)\over |v_s(\rho)|}d\rho}d\eta
\nonumber\\
  &\leq&\displaystyle\int^{R_s}_r|v_s(\eta)|^{-1}e^{-[{\rm Re}
  \lambda-\lambda^*(a)]\int^{\eta}_r|v_s(\rho)|^{-1}d\rho}d\eta
\nonumber\\
  &=&\displaystyle\frac{1}{{\rm Re}\lambda-\lambda^*(a)}\Big[1-e^{-[{\rm Re}
  \lambda-\lambda^*(a)]\int^{R_s}_r|v_s(\rho)|^{-1}d\rho}\Big].
\end{eqnarray*}
  From this estimate and the expression (4.9) we easily see that (4.6) follows.
  This completes the proof of Lemma 4.1. $\quad\Box$
\medskip

  {\bf Corollary 4.2}\ \ {\em The operator $L_0$ generates a positive
  $C_0$-semigroup $e^{tL_0}$ in $C[0,R_s]$ satisfying the following estimate:
  For any $\mu>\lambda_1(a)$ there exists corresponding constant $C_{\mu}>0$ such
  that}
\begin{equation}
  \max_{0\leq r\leq R_s}|e^{tL_0}\phi(r)|\leq C_{\mu}e^{\mu t}
  \max_{0\leq r\leq R_s}|\phi(r)|
  \quad \mbox{for}\;\; \phi\in C[0,R_s],\;\; t\geq 0.
\end{equation}

  {\em Proof}:\ \ By the Hille-Yosida theorem, the assertion (3) of Lemma 4.1
  ensures that $L_0$ generates a $C_0$-semigroup $e^{tL_0}$ in $C[0,R_s]$. The
  assertion (1) of Lemma 4.1 ensures that this semigroup is positive (cf.
  Theorem 1.8 in Chapter VI of \cite{EN}). The assertion (2) of Lemma 4.1 implies
  that the spectral bound of $L_0$ (see Definition 1.12 in Chapter II of \cite{EN}
  for this concept) is not greater than $\lambda_1(a)$: $s(L_0)\leq\lambda_1(a)$. It
  follows by Proposition 1.14 in Chapter VI of \cite{EN} that for any $\mu>
  \lambda_1(a)$ there exists a corresponding constant $C_{\mu}>0$ such that
\begin{equation}
  \max_{0\leq r\leq R_s}|e^{tL_0}1|\leq C_{\mu}e^{\mu t}
  \quad \mbox{for}\;\; t\geq 0.
\end{equation}
  The positivity of $e^{tL_0}$ implies that the comparison principle holds for it.
  Hence, since
$$
  -\max_{0\leq r\leq R_s}|\phi(r)|\leq\phi(r)\leq\max_{0\leq r\leq R_s}|\phi(r)|
$$
  for every $\phi\in C[0,R_s]$, (4.11) is an immediate consequence of (4.12).
  $\quad\Box$
\medskip

  {\bf Lemma 4.3}\ \ {\em The operator $L_0$ also generates a positive
  $C_0$-semigroup $e^{tL_0}$ in $L^1([0,R_s],\\ r^{n-1}dr)$ satisfying the following
  estimate: For any $\phi\in L^1([0,R_s], r^{n-1}dr)$,
\begin{equation}
  \int_0^{R_s}|e^{tL_0}\phi(r)|r^{n-1}dr\leq e^{\lambda_2 t}
  \int_0^{R_s}|\phi(r)|r^{n-1}dr
  \quad \mbox{for}\;\;  t\geq 0,
\end{equation}
  where $\lambda_2=\displaystyle\max_{0\leq r\leq R_s}[g^*(r)+a(r)]$.}
\medskip

  {\em Proof}:\ \ By the density of $C[0,R_s]$ in $L^1([0,R_s], r^{n-1}dr)$, we only
  need to prove the estimate (4.13). Let $\varphi(r,t)=e^{tL_0}\phi(r)$. Then
  $\varphi$ is the solution of the following initial value problem:
\begin{equation}
\left\{
\begin{array}{l}
  \partial_t\varphi=-v_s(r)\partial_r\varphi+a(r)\varphi\;\;
  \mbox{for}\;\; 0<r<R_s,\;\; t>0,\\
  \varphi|_{t=0}=\phi\;\; \mbox{for}\;\; 0\leq r\leq R_s.
\end{array}
\right.
\end{equation}
  Multiplying the first equation in (4.14) with $({\rm sgn}\varphi)r^{n-1}$ and next
  integrating it over $[0,R_s]$, we get
\begin{eqnarray*}
  \displaystyle{d\over dt}\int_0^{R_s}|\varphi(r,t)|r^{n-1}dr
  &=&\displaystyle-\int_0^{R_s}v_s(r)r^{n-1}{\partial\over\partial r}|\varphi(r,t)|dr
  +\int_0^{R_s}a(r)|\varphi(r,t)|r^{n-1}dr
\nonumber\\
  &=&\displaystyle\int_0^{R_s}{\partial\over\partial r}[v_s(r)r^{n-1}]\cdot|\varphi(r,t)|dr
  +\int_0^{R_s}a(r)|\varphi(r,t)|r^{n-1}dr
\nonumber\\
  &=&\displaystyle\int_0^{R_s}[g^*(r)+a(r)]|\varphi(r,t)|r^{n-1}dr
\nonumber\\
  &\leq &\displaystyle\lambda_2\int_0^{R_s}|\varphi(r,t)|r^{n-1}dr.
\end{eqnarray*}
  From this estimate, (4.13) immediately follows. $\quad\Box$
\medskip

  For every nonnegative integer $k$, we let $\hat{\mathscr{L}}_k^+$ be the following
  linear differential-integral operator in $(0,R_s)$: For $\varphi\in C(0,R_s]\cap
  C^1(0,R_s)$,
\begin{equation}
\begin{array}{rcl}
  \hat{\mathscr{L}}_k^+(\varphi)&=&\displaystyle -v_s(r)\varphi'(r)+a_k(r)\varphi(r)
  +{\theta_k}\int_r^{R_s}\!\!\Big(\frac{r}{\rho}\Big)^{n+2(k-1)}
  g_p^*(\rho)p_s'(\rho)\varphi(\rho)d\rho
\\ [0.3cm]
  &&\displaystyle+(1-\theta_k)
  \int_r^{R_s}\!\!g_p^*(\rho)p_s'(\rho)\varphi(\rho)d\rho
   \quad \mbox{for}\;\;0<r<R_s,
\end{array}
\end{equation}
  where
\begin{equation}
  a_k(r)=\frac{k}{r}v_s(r)+g^*(r)-\frac{f_c^*(r)c_s'(r)}{p_s'(r)}.
\end{equation}
  Since $g^*(0)<0$, $v_s'(0)=\displaystyle\frac{1}{n}g^*(0)<0$ and $\displaystyle
  \frac{f_c^*(r)c_s'(r)}{p_s'(r)}
  >0$, we easily see that $a_k(0)=\displaystyle\lim_{r\to 0^+}a_k(r)<0$ (note that
  since $p_s''(0)\neq 0$ in case $p_s'(0)=0$ (see (4.10) in \cite{CuiFri1}), this
  limit exists). We have also that $a_k(R_s)<0$. Indeed, since $f_c^*(r)c_s'(r)+
  f_p^*(r)p_s'(r)=(v_s(r)p_s'(r))'=v_s'(r)p_s'(r)+v_s(r)p_s''(r)$, $v_s(R_s)=0$ and
  $v_s'(R_s)=g^*(R_s)$, we see that
\begin{equation*}
  a_k(R_s)=g^*(R_s)-v_s'(R_s)+f_p^*(R_s)=f_p^*(R_s)<0.
\end{equation*}
  We denote
\begin{equation*}
  \mu_k=\max\{a_k(0),a_k(R_s)\}.
\end{equation*}
  The above argument shows that $\mu_k<0$, $k=0,1,2,\cdots$, and, in fact,
\begin{equation*}
  \mu_k\leq\max\{\Big(1+\frac{k}{n}\Big)g^*(0),f_p^*(R_s)\}
  \leq\max\{g^*(0),f_p^*(R_s)\}\equiv\mu_0^*<0, \quad
  k=0,1,2,\cdots.
\end{equation*}

  {\bf Lemma 4.4}\ \ {\em Let ${\rm Re}\lambda>\mu_k$. For any $h\in C[0,R_s]$ the
  equation
\begin{equation}
  \lambda\varphi-\hat{\mathscr{L}}_k^+(\varphi)=h \quad \mbox{in}\;\;(0,R_s)
\end{equation}
  has a unique solution $\varphi=\varphi_{\lambda}\in C[0,R_s]\cap C^1(0,R_s)$, and
  there exists a constant $C_{k,\lambda}>0$ such that
\begin{equation}
  \max_{0\leq r\leq R_s}|\varphi_{\lambda}(r)|\leq C_{k,\lambda}\max_{0\leq r\leq R_s}|h(r)|.
\end{equation}
  Moreover, if $\lambda$ is real and $h\geq 0$ then also $\varphi_{\lambda}\geq 0$.}
\medskip

  {\em Proof}:\ \ We fulfill the proof through three steps.

  {\em Step 1}:\ \ We first prove that the equation (4.11) has a unique solution in
  the class $C(0,R_s]\cap C^1(0,R_s)$. To this end, we explicitly write out the equation
  (4.11) as follows:
\begin{eqnarray}
  \displaystyle v_s(r)\varphi'(r)&+&[\lambda-a_k(r)]\varphi(r)
  -{\theta_k}\int_r^{R_s}\!\!\Big(\frac{r}{\rho}\Big)^{n+2(k-1)}
  g_p^*(\rho)p_s'(\rho)\varphi(\rho)d\rho
\nonumber\\
  &-&\displaystyle(1-\theta_k)
  \int_r^{R_s}\!\!g_p^*(\rho)p_s'(\rho)\varphi(\rho)d\rho=h(r).
\end{eqnarray}
  Let $W_{\lambda}(r)$ be as before but with $a(r)$ replaced with $a_k(r)$.
  By rewriting the above equation in the form
$$
\begin{array}{rcl}
  \displaystyle\frac{d}{dr}\Big(W_{\lambda}(r)\varphi(r)\Big)
  &=&\displaystyle{W_{\lambda}(r)\over v_s(r)}\Big[h(r)+{\theta_k}
  \int_r^{R_s}\!\!\Big(\frac{r}{\rho}\Big)^{n+2(k-1)}g_p^*(\rho)p_s'(\rho)\varphi(\rho)d\rho
\\ [0.3cm]
  &&+\displaystyle(1-\theta_k)
  \int_r^{R_s}\!\!g_p^*(\rho)p_s'(\rho)\varphi(\rho)d\rho\Big],
\end{array}
$$
  we easily see that, as far as solutions which are bounded near $r=R_s$ are concerned,
  the differential-integral equation (4.22) is equivalent to the following integral
  equation:
\begin{equation}
\begin{array}{rcl}
  \displaystyle\varphi(r)&=&\displaystyle{1\over W_{\lambda}(r)}
  \int^{R_s}_r{W_{\lambda}(\eta)\over |v_s(\eta)|}
  \Big[h(\eta)+{\theta_k}\int_{\eta}^{R_s}\!\!
  \Big(\frac{\eta}{\rho}\Big)^{n+2(k-1)}g_p^*(\rho)p_s'(\rho)\varphi(\rho)d\rho
\\ [0.3cm]
  &&+\displaystyle(1-\theta_k)
  \int_{\eta}^{R_s}\!\!g_p^*(\rho)p_s'(\rho)\varphi(\rho)d\rho\Big]d\eta.
\end{array}
\end{equation}
  It follows from a standard contraction mapping argument (similar to that used in the
  proof of Theorem 5.3 (1) of \cite{ChenCuiF}) that there exists a sufficiently small
  $\delta>0$ such that (4.19) has a unique bounded solution in the interval
  $(R_s-\delta,R_s)$, such that $\varphi\in C(R_s-\delta,R_s]\cap C^1(R_s-\delta,R_s)$,
  and
\begin{equation}
  \varphi(r)=\frac{h(R_s)}{v_s'(R_s)}\frac{1}{\alpha_1}[1+o(1)]
  =\frac{h(R_s)[1+o(1)]}{\lambda-a_k(R_s)} \quad \mbox{as}\;\; r\to R_s^-.
\end{equation}
  Since $v_s(r)\neq 0$ for $0<r<R_s$, the equation (4.19) is a regular linear
  differential-integral equation at any point in $(0,R_s)$, so that by standard ODE
  theory we can uniquely extend the solution to the whole interval $(0,R_s]$ such that
  $\varphi\in C(0,R_s]\cap C^1(0,R_s)$. This fulfills the task of the first step.

  We note that if $\lambda$ is real, $\lambda>\mu_k$ and $h\geq 0$, then also
  $\varphi\geq 0$. Indeed, since $\lambda$ is real, we have $W_{\lambda}(r)>0$ for
  $0<r<R_s$. If $h(R_s)>0$ then by (4.21) we see that $\varphi(R_s)>0$. Let $r_0$
  be the smallest number such that $\varphi(r)>0$ for $r_0<r\leq R_s$. Then by (4.20)
  we must have $r_0=0$. The assertion for the case $h(R_s)=0$ follows from a limit
  argument.

  {\em Step 2}:\ \ We next prove that the solution ensured by the above step satisfies
$$
  \int_0^{R_s}\!\!g_p^*(\rho)p_s'(\rho)|\varphi(\rho)|d\rho<\infty.
$$
  To prove this assertion we note that from (4.20) we have
$$
\begin{array}{rl}
  \displaystyle |\varphi(r)|\leq &
  \displaystyle{1\over |W_{\lambda}(r)|}\int^{R_s}_r{|W_{\lambda}(\eta)|
  \over |v_s(\eta)|}\Big[|h(\eta)|+{\theta_k}
  \int_{\eta}^{R_s}\Big(\frac{\eta}{\rho}\Big)^{n+2(k-1)}g_p^*(\rho)p_s'(\rho)|\varphi(\rho)|d\rho
\\ [0.3cm]
  &\displaystyle +(1-\theta_k)
  \int_\eta^{R_s}g_p^*(\rho)p_s'(\rho)|\varphi(\rho)|d\rho\Big]d\eta
\\ [0.3cm]
  \leq &\displaystyle{1\over |W_{\lambda}(r)|}
  \int^{R_s}_r{|W_{\lambda}(\eta)|\over |v_s(\eta)|}\Big[|h(\eta)|
  +\int_\eta^{R_s}g_p^*(\rho)p_s'(\rho)|\varphi(\rho)|d\rho\Big]d\eta.
\end{array}
$$
  It follows that for any $0<r<r'\leq R_s$ we have
\begin{eqnarray*}
  \displaystyle\int^{r'}_rg_p^*(\xi)p_s'(\xi)|\varphi(\xi)|d\xi &\leq &
  \displaystyle\int^{r'}_r\!\!\int^{R_s}_{\xi}
  {g_p^*(\xi)p_s'(\xi)|W_{\lambda}(\eta)||h(\eta)|\over
  |W_{\lambda}(\xi)||v_s(\eta)|}d\eta d\xi
\\
  &&\displaystyle+C\int^{r'}_r\!\!\int^{R_s}_{\xi}\!\!\int_{\eta}^{R_s}
  {g_p^*(\xi)p_s'(\xi)|W_{\lambda}(\eta)|\over|W_{\lambda}(\xi)||v_s(\eta)|}
  g_p^*(\rho)p_s'(\rho)|\varphi(\rho)|d\rho d\eta d\xi
\\
  &\leq &\displaystyle C\Big(\int^{r'}_r\!\!\int^{R_s}_{\xi}
  {p_s'(\xi)|W_{\lambda}(\eta)|\over |W_{\lambda}(\xi)||v_s(\eta)|}
  d\eta d\xi\Big)\max_{0\leq r\leq R_s}|h(r)|
\\
  &&\displaystyle+C\Big(\int^{r'}_r\!\!\int^{R_s}_{\xi}
  {p_s'(\xi)|W_{\lambda}(\eta)|\over |W_{\lambda}(\xi)||v_s(\eta)|}
  d\eta d\xi\Big)\Big(\int_r^{R_s}g_p^*(\rho)p_s'(\rho)|\varphi(\rho)|d\rho\Big).
\end{eqnarray*}
  From (4.4) and (4.5) we have
\begin{equation}
  C_1r^{-{\rm Re}\alpha_0}(R_s-r)^{{\rm Re}\alpha_1}\leq |W_{\lambda}(r)|
  \leq C_2r^{-{\rm Re}\alpha_0}(R_s-r)^{{\rm Re}\alpha_1} \quad \mbox{for}\;\; 0<r<R_s,
\end{equation}
  where $0<C_1<C_2$, ${\rm Re}\alpha_0>0$ and ${\rm Re}\alpha_1>0$. Using these
  estimates and (3.2) we have
\begin{eqnarray*}
  \displaystyle\int^{R_s}_0\!\!\int^{R_s}_{\xi}
  {p_s'(\xi)|W_{\lambda}(\eta)|\over |W_{\lambda}(\xi)||v_s(\eta)|}
  d\eta d\xi
  &=&\displaystyle\Big(\int^{{R_s\over 2}}_0\!\!\int^{{R_s\over 2}}_{\xi}
  +\int^{{R_s\over 2}}_0\!\!\int^{R_s}_{{R_s\over 2}}
  +\int^{R_s}_{{R_s\over 2}}\!\!\int^{R_s}_{\xi}\Big)
  {p_s'(\xi)|W_{\lambda}(\eta)|\over |W_{\lambda}(\xi)||v_s(\eta)|}
  d\eta d\xi
\\
  &\leq &\displaystyle C\int^{{R_s\over 2}}_0\!\!\int^{{R_s\over 2}}_{\xi}
  p_s'(\xi)\xi^{{\rm Re}\alpha_0}\eta^{-{\rm Re}\alpha_0-1}d\eta d\xi
\\
  &&\displaystyle+C\int^{{R_s\over 2}}_0\!\!\int^{R_s}_{{R_s\over 2}}
  p_s'(\xi)\xi^{{\rm Re}\alpha_0}(R_s-\eta)^{{\rm Re}\alpha_1-1}d\eta d\xi
\\
  &&\displaystyle+C\int^{R_s}_{{R_s\over 2}}\!\!\int^{R_s}_{\xi}
  (R_s-\xi)^{-{\rm Re}\alpha_1}(R_s-\eta)^{{\rm Re}\alpha_1-1}d\eta d\xi
\\
  &&<\infty.
\end{eqnarray*}
  Hence there exists a constant $\delta>0$ independent of $k$ such that if $0<r'-r\leq\delta$
  then
$$
  C\int^{r'}_r\!\!\int^{R_s}_{\xi}
  {p_s'(\xi)|W_{\lambda}(\eta)|\over |W_{\lambda}(\xi)||v_s(\eta)|}
  d\eta d\xi\leq\frac{1}{2},
$$
  which implies that
$$
  \int^{r'}_rg_p^*(\xi)p_s'(\xi)|\varphi(\xi)|d\xi
  \leq\max_{0\leq r\leq R_s}|h(r)|
  +\int_{r'}^{R_s}g_p^*(\xi)p_s'(\xi)|\varphi(\xi)|d\xi.
$$
  Hence, by dividing the interval $[0,R_s]$ into finite number (depending on $k$
  because $W_{\lambda}$ depends on $k$) of subintervals and using an iteration
  argument, we see that there exists a constant $C_{k,\lambda}>0$ depending on $k$
  such that
$$
  \int^{R_s}_0g_p^*(\rho)p_s'(\rho)|\varphi(\rho)|d\rho
  \leq C_{k,\lambda}\max_{0\leq r\leq R_s}|h(r)|.
$$
  This fulfills the task of the second step.

  {\em Step 3}:\ \ From the assertion obtained in the above step, it follows that
  the function
$$
  h_1(r)={\theta_k}
  \int_r^{R_s}\!\!\Big(\frac{r}{\rho}\Big)^{n+2(k-1)}g_p^*(\rho)p_s'(\rho)\varphi(\rho)d\rho
  +(1-\theta_k)
  \int_r^{R_s}\!\!g_p^*(\rho)p_s'(\rho)\varphi(\rho)d\rho
$$
  belongs to $C[0,R_s]$, and $\displaystyle\max_{0\leq r\leq R_s}|h_1(r)|\leq
  C_{k,\lambda}\max_{0\leq r\leq R_s}|h(r)|$. Hence, by rewritting the equation
  (4.19) into the form
$$
  v_s(r)\varphi'(r)+[\lambda-a_k(r)]\varphi(r)=h(r)+h_1(r)
$$
  and applying the assertion (2) of Lemma 4.1, we get the desired assertion.
  This completes the proof of Lemma 4.3. $\quad\Box$
\medskip

  {\bf Corollary 4.5}\ \ {\em The operator $\hat{\mathscr{L}}_k^+$ generates a positive
  $C_0$-semigroup $e^{t\hat{\mathscr{L}}_k^+}$ in $C[0,R_s]$ satisfying the following
  estimate: For any $\mu>\mu_k$ there exists corresponding constant $C_{k,\mu}>0$
  such that}
\begin{equation}
  \max_{0\leq r\leq R_s}|e^{t\hat{\mathscr{L}}_k^+}\phi(r)|\leq C_{k,\mu}e^{\mu t}
  \max_{0\leq r\leq R_s}|\phi(r)|
  \quad \mbox{for}\;\; \phi\in C[0,R_s],\;\; t\geq 0.
\end{equation}

  {\em Proof}:\ \ We note that $\hat{\mathscr{L}}_k^+=\hat{\mathscr{L}}_k^0+\mathscr{B}_k$,
  where
\begin{equation}
  \hat{\mathscr{L}}_k^0(\phi)=-v_s(r)\phi'(r)+a_k(r)\phi(r)
  \quad \mbox{for}\;\; \phi\in C[0,R_s],
\end{equation}
  and $\mathscr{B}_k$ is the integral part of $\hat{\mathscr{L}}_k^+$.
  By Corollary 4.2 we see that $\hat{\mathscr{L}}_k^0$ generates a positive
  $C_0$-semigroup in $C[0,R_s]$. Since clearly $\mathscr{B}_k$ is a positive bounded
  linear operator in $C[0,R_s]$, by a standard perturbation theorem (cf. Corollary
  1.11 in Chapter VI of \cite{EN}) for $C_0$-semigroups we see that
  $\hat{\mathscr{L}}_k^+$ also generates a positive $C_0$-semigroup in $C[0,R_s]$.
  By Lemma 4.4 we see that the spectral bound of $\hat{\mathscr{L}}_k^+$ is
  not greater than $\mu_k$: $s(\hat{\mathscr{L}}_k^+)\leq\mu_k$. Hence by a
  similar argument as in the proof of Corollary 4.2 we obtain the estimate
  (4.23). $\quad\Box$
\medskip

  Let $J$ be the following continuous linear functional in $L^1[0,R_s]$:
$$
  J(\phi)=\int_0^{R_s}g_p^*(r)p_s'(r)\phi(r)dr \quad \mbox{for} \quad
  \phi\in L^1[0,R_s].
$$

  {\bf Lemma 4.6}\ \ {\em Assume that the conditions in $(1.25)$ are satisfied. For $\psi_0\in C[0,R_s]$
  we let $\psi=e^{t\hat{\mathscr{L}}_k^+}\psi_0$. There exists a constant $\kappa_0>0$
  independent of $k$ such that the following assertion holds for any $\psi_0\in
  C[0,R_s]$: If $\psi_0\geq 0$ then
\begin{equation}
  \frac{d}{dt}J(\psi(\cdot,t))\leq -\kappa_0 J(\psi) \quad \mbox{for}\;\; t\geq 0,
\end{equation}
  and for general $\psi_0\in C[0,R_s]$ we have}
\begin{equation}
  \frac{d}{dt}J(|\psi(\cdot,t)|)\leq -\kappa_0 J(|\psi|) \quad \mbox{for}\;\; t\geq 0.
\end{equation}

  {\em Proof}:\ \ Since the semigroup $e^{t\hat{\mathscr{L}}_k^+}$ is positive, we
  see that $\psi_0\geq 0$ implies $\psi\geq 0$. To prove (4.25) we note that $\psi$
  is a solution of the following equation:
\begin{equation}
  \partial_t\psi=\hat{\mathscr{L}}_k^+\psi\;\;
  \mbox{for}\;\; 0<r<R_s,\;\; t>0.
\end{equation}
  Using this fact we compute
\begin{eqnarray}
  \frac{d}{dt}J(\psi)&=&\int_0^{R_s}g_p^*(r)p_s'(r)\partial_t\psi(r,t)dr
\nonumber\\
  &=&\displaystyle -\int_0^{R_s}v_s(r)g_p^*(r)p_s'(r)\partial_r\psi(r,t)dr
  +\int_0^{R_s}a_k(r)g_p^*(r)p_s'(r)\psi(r,t)dr
\nonumber\\
  &&+{\theta_k}\int_0^{R_s}\!\!\int_r^{R_s}\!\!\Big(\frac{r}{\rho}\Big)^{n+2(k-1)}
  g_p^*(r)p_s'(r)g_p^*(\rho)p_s'(\rho)\varphi(\rho)d\rho dr
\nonumber\\
  &&\displaystyle+(1-\theta_k)
  \int_0^{R_s}\!\!\int_r^{R_s}\!g_p^*(r)p_s'(r)g_p^*(\rho)p_s'(\rho)\varphi(\rho)d\rho dr
\nonumber\\
  &=&\displaystyle\int_0^{R_s}v_s(r)g_p^{*'}(r)p_s'(r)\psi(r,t)dr
  +\int_0^{R_s}g_p^*(r)[v_s(r)p_s'(r)]'\psi(r,t)dr
\nonumber\\
  && +\int_0^{R_s}\Big\{a_k(\rho)+{\theta_k}
  \int_0^{\rho}\!\!\Big(\frac{r}{\rho}\Big)^{n+2(k-1)}g_p^*(r)p_s'(r)dr
\nonumber\\
  &&\displaystyle+(1-\theta_k)
  \int_0^{\rho}g_p^*(r)p_s'(r)dr\Big\}g_p^*(\rho)p_s'(\rho)\psi(\rho,t)d\rho
\nonumber\\
  &\leq &\displaystyle\int_0^{R_s}\tilde{a}_k(r)g_p^*(r)p_s'(r)\psi(r,t)dr,
\end{eqnarray}
  where
$$
  \tilde{a}_k(r)=\frac{k}{r}v_s(r)+\frac{g_p^{*'}(r)}{g_p^*(r)}v_s(r)
  +f_p^*(r)+g^*(r)+\int_0^rg_p^*(\rho)p_s'(\rho)d\rho.
$$
  It is clear that the first two terms on the right-hand side of the above equality
  are negative for $0<r<R_s$. In what follows we prove that if the conditions in $(1.25)$ are satisfied
  then
\begin{equation}
  f_p^*(r)+g^*(r)+\int_0^rg_p^*(\rho)p_s'(\rho)d\rho<0 \quad
  \mbox{for}\;\; 0\leq r\leq R_s.
\end{equation}
  We first note that, since $g_p^*(r)=K_M(c_s(r))$ is monotone increasing in $r$, we
  have
\begin{eqnarray*}
  \int_0^rg_p^*(\rho)p_s'(\rho)d\rho &\leq & g_p^*(r)\int_0^rp_s'(\rho)d\rho
  =g_p^*(r)[p_s(r)-p_s(0)]
\nonumber\\
  &=&g^*(r)+K_D(c_s(r))-p_s(0)K_M(c_s(r)).
\end{eqnarray*}
  Hence
\begin{eqnarray*}
  f_p^*(r)+g^*(r)&+&\int_0^rg_p^*(\rho)p_s'(\rho)d\rho
  \leq f_p^*(r)+2g^*(r)+K_D(c_s(r))-p_s(0)K_M(c_s(r))
\nonumber\\
  &=&K_M(c_s(r))-K_N(c_s(r))-K_D(c_s(r))-p_s(0)K_M(c_s(r))
\nonumber\\
  &=&K_M(c_s(r))[1-p_s(0)]-K_N(c_s(r))-K_D(c_s(r)).
\end{eqnarray*}
  Note that the conditions in (1.25) imply that the function $c\mapsto
  \displaystyle\frac{K_N(c)}{K_M(c)}=\frac{(k_P-k_Q)c+k_Q}{(k_B-k_D)c+k_D}$ is monotone
  increasing for $c>0$. Hence
\begin{eqnarray*}
  p_s(0)&=&\frac{1}{2K_M(c_s(0))}\{[K_M(c_s(0))-K_N(c_s(0))]
  +\sqrt{[K_M(c_s(0))-K_N(c_s(0))]^2+4K_M(c_s(0))K_P(c_s(0))}\}
\nonumber\\
  &>&1-\frac{K_N(c_s(0))}{K_M(c_s(0))}\geq 1-\frac{K_N(c_s(r))}{K_M(c_s(r))}.
\end{eqnarray*}
  From these estimates we see that (4.29) follows. Having proved (4.29), we see
  that
$$
  \kappa_0=-\max_{0\leq r\leq R_s}\Big(f_p^*(r)+g^*(r)
  +\int_0^rg_p^*(\rho)p_s'(\rho)d\rho\Big)>0,
$$
  and
$$
  \tilde{a}_k(r)\leq -\kappa_0 \quad \mbox{for} \;\; 0\leq r\leq R_s.
$$
  Using this result and (4.28), we see that (4.25) follows. To prove (4.26) we
  multiply the equation in (4.27) with ${\rm sgn}\psi$, which yields the following
  relation:
$$
  \partial_t|\psi|\leq\hat{\mathscr{L}}_k^+|\psi|\;\;
  \mbox{for}\;\; 0<r<R_s,\;\; t>0.
$$
  Using this fact and a similar argument as above we obtain (4.26). This completes
  the proof of Lemma 4.6. $\quad\Box$

\section{Decay estimates for the equation
  $\displaystyle\partial_t\varphi=\tilde{\mathscr{L}}_k(\varphi)$ for large $k$}
\setcounter{equation}{0}

\hskip 2em
  In this and the next sections we establish decay estimates for the solution of
  the following initial value problem:
\begin{equation}
\left\{
\begin{array}{l}
  \partial_t\varphi=\tilde{\mathscr{L}}_k(\varphi)\;\;
  \mbox{for}\;\; 0<r<R_s,\;\; t>0,\\
  \varphi|_{t=0}=\varphi_0\;\; \mbox{for}\;\; 0\leq r\leq R_s.
\end{array}
\right.
\end{equation}
  In what follows we consider the case that $k$ is sufficiently large; the rest cases
  will be treated in the next subsection. We denote
$$
  \nu_0=\max\{f_p^*(0),f_p^*(R_s)\}.
$$
  Since $f_p^*(r)<0$ for all $0\leq r\leq R_s$, we see that $\nu_0<0$.

\medskip
  {\bf Lemma 5.1}\ \ {\em For any $\mu>\nu_0$ there exist corresponding positive
  integer $k_{\mu}$ and positive constant $C=C_{\mu}$ such that for any $k\geq k_{\mu}$
  and $\varphi_0\in C[0,R_s]$, the solution of the initial value problem $(5.1)$
  satisfies the following estimate:}
\begin{equation}
  \max_{0\leq r\leq R_s}|\varphi(r,t)|
  \leq C\max_{0\leq r\leq R_s}|\varphi_0(r)|e^{\mu t} \quad \mbox{for}\;\; t\geq 0.
\end{equation}

  {\em Proof}:\ \ Let $\mathscr{L}_0$ be the following unbounded closed linear operator
  in $C[0,R_s]$ with domain $D(\mathscr{L}_0)=\{\varphi\in C[0,R_s]\cap C^1(0,R_s):
  v_s(r)\varphi'(r)\in C[0,R_s]\}$: For $\phi\in D(\mathscr{L}_0)$,
$$
  \mathscr{L}_0(\phi)=-v_s(r)\phi'(r)+f_p^*(r)\phi(r) \quad
  \mbox{for}\;\; 0\leq r\leq R_s.
$$
  By Corollary 4.2 we see that $\mathscr{L}_0$ generates a positive $C_0$-semigroup
  $e^{t\mathscr{L}_0}$ in $C[0,R_s]$ satisfying the following estimate: For any
  $\mu>\nu_0$,
\begin{equation}
  \|e^{t\mathscr{L}_0}\|_{L(C[0,R_s])}\leq C_{\mu}e^{\mu t}
   \quad \mbox{for}\;\; t\geq 0.
\end{equation}
  Now for each integer $k\geq 2$ we denote by $\mathscr{K}_k$ the following bounded
  linear operator in $C[0,R_s]$: For any $\varphi\in C[0,R_s]$,
\begin{eqnarray}
  \qquad\mathscr{K}_k(\varphi)(r)&=&\displaystyle
  r^{k-1}p_s'(r)\Big[{\theta_k}
  \int_r^{R_s}\rho^{-k+1}g_p^*(\rho)\varphi(\rho)d\rho
  -\frac{1-\theta_k}{r^{n+2(k-1)}}
  \int_0^r\rho^{n+k-1}g_p^*(\rho)\varphi(\rho)d\rho
\nonumber\\
  &&\displaystyle+\frac{1-\theta_k}{R_s^{n+2(k-1)}}
  \int_0^{R_s}\rho^{n+k-1}g_p^*(\rho)\varphi(\rho)d\rho\Big].
\end{eqnarray}
  It is easy to see that there exists a positive constant $C$ independent of $k$
  such that for any $k\geq 3$,
$$
  \max_{0\leq r\leq R_s}|\mathscr{K}_k(\varphi)(r)|\leq
  Ck^{-1}\max_{0\leq r\leq R_s}|\varphi(r)| \quad
  \mbox{for any}\;\; \varphi\in C[0,R_s].
$$
  Since $\tilde{\mathscr{L}}_k=\mathscr{L}_0+\mathscr{K}_k$ ($k=2,3,\cdots$), by
  using a standard perturbation theorem for $C_0$-semigroups we deduce from the
  above estimate and (5.3) that for any $\mu>\nu_0$,
\begin{equation*}
  \|e^{t\tilde{\mathscr{L}}_k}\|_{L(C[0,R_s])}\leq C_{\mu}e^{(\mu+C_{\mu}Ck^{-1})t}
   \quad \mbox{for}\;\; t\geq 0.
\end{equation*}
  The desired assertion immediately follows from this result. This proves
  Lemma 5.1. $\quad\Box$
\medskip

  For every $\alpha\geq 1$ we denote
$$
  \mu_{\alpha}^*=\max_{0\leq r\leq R_s}\Big(f_p^*(r)+\frac{1}{\alpha}g^*(r)\Big).
$$
  Note that
\begin{eqnarray*}
  f_p^*(r)+\frac{1}{\alpha}g^*(r)&=&[K_M(c_s(r))-K_N(c_s(r))]-2K_M(c_s(r))p_s(r)
  +\frac{1}{\alpha}[K_M(c_s(r))p_s(r)-K_D(c_s(r))]
\nonumber\\
  &<&-K_M(c_s(r))p_s(r)+\frac{1}{\alpha}[K_M(c_s(r))p_s(r)-K_D(c_s(r))]
\nonumber\\
  &=&-\Big(1-\frac{1}{\alpha}\Big)K_M(c_s(r))p_s(r)-\frac{1}{\alpha}K_D(c_s(r)).
\nonumber\\
  &\leq& 0 \quad \mbox{for}\;\; 0\leq r\leq R_s.
\end{eqnarray*}
  Hence $\mu_{\alpha}^*<0$ for all $\alpha\geq 1$.
\medskip

  {\bf Lemma 5.2}\ \ {\em For any $\alpha\geq 1$ and $\mu>\mu_{\alpha}^*$ there exist
  corresponding positive integer $k_{\alpha,\mu}$ and positive constant $C=C_{\alpha,\mu}$
  such that for any $k\geq k_{\mu}$ and $\varphi_0\in L^{\alpha}([0,R_s];r^{n-1}dr)$,
  the solution of the initial value problem $(5.1)$ satisfies the following estimate:}
\begin{equation}
  \Big(\int_0^{R_s}|\varphi(r,t)|^{\alpha}r^{n-1}dr\Big)^{\frac{1}{\alpha}}
  \leq C\Big(\int_0^{R_s}|\varphi_0(r)|^{\alpha}r^{n-1}dr\Big)^{\frac{1}{\alpha}}e^{\mu t}
  \quad \mbox{for}\;\; t\geq 0.
\end{equation}

  {\em Proof}:\ \ We first establish a $L^{\alpha}$-estimate for $e^{t\mathscr{L}_0}$.
  Given $\varphi_0\in C[0,R_s]$ we let $\varphi(r,t)=e^{t\mathscr{L}_0}\varphi_0(r)$.
  Then $\varphi$ is a solution of the following initial value problem:
\begin{equation}
\left\{
\begin{array}{l}
  \partial_t\varphi=\mathscr{L}_0(\varphi)\;\;
  \mbox{for}\;\; 0<r<R_s,\;\; t>0,\\
  \varphi|_{t=0}=\varphi_0\;\; \mbox{for}\;\; 0\leq r\leq R_s.
\end{array}
\right.
\end{equation}
  By a standard argument we have
\begin{eqnarray*}
  \frac{1}{\alpha}\frac{d}{dt}\int_0^{R_s}r^{n-1}|\varphi(r,t)|^{\alpha}dr
  &=&-\frac{1}{\alpha}\int_0^{R_s}r^{n-1}v_s(r)
  \frac{\partial}{\partial r}|\varphi(r,t)|^{\alpha}dr
  +\int_0^{R_s}r^{n-1}f_p^*(r)|\varphi(r,t)|^{\alpha}dr
\nonumber\\
  &=&\frac{1}{\alpha}\int_0^{R_s}\frac{\partial}{\partial r}\Big(r^{n-1}v_s(r)\Big)
  |\varphi(r,t)|^{\alpha}dr
  +\int_0^{R_s}r^{n-1}f_p^*(r)|\varphi(r,t)|^{\alpha}dr
\nonumber\\
  &=&\int_0^{R_s}\Big(\frac{1}{\alpha}g^*(r)+f_p^*(r)\Big)r^{n-1}|\varphi(r,t)|^{\alpha}dr
\nonumber\\
  &\leq&\mu_{\alpha}^*\int_0^{R_s}r^{n-1}|\varphi(r,t)|^{\alpha}dr.
\end{eqnarray*}
  Hence
\begin{equation*}
  \Big(\int_0^{R_s}|\varphi(r,t)|^{\alpha}r^{n-1}dr\Big)^{\frac{1}{\alpha}}
  \leq\Big(\int_0^{R_s}|\varphi_0(r)|^{\alpha}r^{n-1}dr\Big)^{\frac{1}{\alpha}}e^{\mu_{\alpha}^* t}
  \quad \mbox{for}\;\; t\geq 0.
\end{equation*}
  This means that, by using the abbreviation $L^{\alpha}$ for $L^{\alpha}([0,R_s];r^{n-1}dr)$,
  we have
\begin{equation}
  \|e^{t\mathscr{L}_0}\varphi_0\|_{L^{\alpha}}
  \leq\|\varphi_0\|_{L^{\alpha}}e^{\mu_{\alpha}^* t}
  \quad \mbox{for}\;\; t\geq 0.
\end{equation}
  Next, it is not hard to prove that there exists a constant $C_{\alpha}>0$ independent
  of $k$ such that for any $k\geq 3$ and $\phi\in L^{\alpha}([0,R_s];r^{n-1}dr)$,
\begin{equation}
  \|\mathscr{K}_k\phi\|_{L^{\alpha}}
  \leq C_{\alpha}k^{-1}\|\phi\|_{L^{\alpha}}.
\end{equation}
  Since $\tilde{\mathscr{L}}_k=\mathscr{L}_0+\mathscr{K}_k$ ($k=2,3,\cdots$), from
  (5.7), (5.8) and a standard perturbation theorem for $C_0$-semigroups we
  see that the desired assertion follows. $\quad\Box$

\section{Decay estimates for the equation
  $\displaystyle\partial_t\varphi=\tilde{\mathscr{L}}_k(\varphi)$ for small $k$}
\setcounter{equation}{0}

\hskip 2em
  A similar estimate for small $k$ is much more involved. In what follows we
  consider this case.

  For every nonnegative integer $k$, we denote by $\hat{\mathscr{L}}_k$ the following
  differential-integral operator:
\begin{eqnarray}
  \hat{\mathscr{L}}_k\phi(r)&=&\displaystyle-v_s(r)\phi'(r)+a_k(r)\phi(r)
  +{\theta_k}
  \int_r^{R_s}\Big(\frac{r}{\rho}\Big)^{n+2(k-1)}g_p^*(\rho)p_s'(\rho)\phi(\rho)d\rho
\nonumber\\
  &&\displaystyle+(1-\theta_k)\Big(\frac{r}{R_s}\Big)^{n+2(k-1)}
  \int_0^{R_s}g_p^*(\rho)p_s'(\rho)\phi(\rho)d\rho
\nonumber\\
  &&\displaystyle-(1-\theta_k)\int_0^rg_p^*(\rho)p_s'(\rho)\phi(\rho)d\rho,
\end{eqnarray}
  where $a_k(r)$ is as before (see (4.16)). By using the relations
$$
  v_s(r)p_s'(r)=f(c_s(r),p_s(r)),
$$
$$
  v_s'(r)+\frac{n-1}{r}v_s(r)=g^*(r),
$$
  we can easily get the following relation:
\begin{equation}
  \tilde{\mathscr{L}}_k[r^{-(n+k-1)}p_s'(r)\phi(r)]
  =r^{-(n+k-1)}p_s'(r)\hat{\mathscr{L}}_k\phi(r).
\end{equation}
  It follows that if we let $\varphi(r,t)=r^{-(n+k-1)}p_s'(r)\psi(r,t)$ and
  $\varphi_0(r)=r^{-(n+k-1)}p_s'(r)\psi_0(r)$, then $\varphi$ is a solution of the
  problem (5.1) if and only if $\psi$ is a solution of the following problem:
\begin{equation}
\left\{
\begin{array}{l}
  \partial_t\psi=\hat{\mathscr{L}}_k\psi\;\;
  \mbox{for}\;\; 0<r<R_s,\;\; t>0,\\
  \psi|_{t=0}=\psi_0.
\end{array}
\right.
\end{equation}
  Note that by denoting
$$
  e_k(r)=(1-\theta_k)
  \Big[1-\Big(\frac{r}{R_s}\Big)^{n+2(k-1)}\Big],
$$
  we have
$$
  \hat{\mathscr{L}}_k\phi(r)=\hat{\mathscr{L}}_k^+\phi(r)-J(\phi)e_k(r)
  \quad \mbox{for}\;\; \phi\in D(\hat{\mathscr{L}}_k)=D(\hat{\mathscr{L}}_k^+).
$$

  {\bf Lemma 6.1}\ \ {\em Let $\psi(r,t)$ be the solution of the problem $(6.3)$,
  $\kappa(r,t)=e^{t\hat{\mathscr{L}}_k^+}e_k(r)$, and $\tilde{\psi}(r,t)=
  e^{t\hat{\mathscr{L}}_k^+}\psi_0(r)$. Let $\Psi(t)=J(\psi)$, $K(t)=J(\kappa)$
  and $\tilde{\Psi}(t)=J(\tilde{\psi})$. Then the following relation holds:}
\begin{equation}
  \Psi(t)+\int^t_0\Psi(\tau)K(t-\tau)d\tau=\tilde{\Psi}(t) \quad \mbox{for}\;\; t\geq 0.
\end{equation}

  {\em Proof}:\ \ The proof is similar to that of Lemma 8.2 of \cite{ChenCuiF},
  so that is omitted. $\quad\Box$
\medskip

  {\bf Lemma 6.2}\ \ {\em Let $K\in C^1[0,\infty)$ and assume that
$$
   K(t)\geq 0, \quad {d\over dt}\big(e^{\sigma t}K(t)\big)\leq 0 \;\;
   {\it for}\;\;t\geq 0
$$
  for some real constant $\sigma$. Then for any $\tilde{\Psi}\in C[0,\infty)$, the unique
  solution $\Psi$ of the Volterra integral equation $(6.4)$ satisfies the following
  estimate:}
\begin{equation*}
   |\Psi(t)-\tilde{\Psi}(t)|\leq K(0)\int^t_0e^{-\sigma(t-\tau)}
   |\tilde{\Psi}(\tau)|d\tau\;\;\;{\it for}\;\;t\geq 0.
\end{equation*}

  {\em Proof}:\ \ See Lemma 8.3 of \cite{ChenCuiF}. $\quad\Box$
\medskip

  {\bf Lemma 6.3}\ \ {\em Assume that the conditions in $(1.25)$ are satisfied. There exists a constant
  $\mu^*<0$ independent of $k$ such that for every nonnegative integer $k$ and any
  $\mu>\mu^*$ there exists corresponding constant $C=C_{k,\mu}>0$ such that for the
  solution of the problem $(6.3)$ the following estimate holds:}
\begin{equation}
  \max_{0\leq r\leq R_s}|\psi(r,t)|\leq Ce^{\mu t}\max_{0\leq r\leq R_s}|\psi_0(r)|
  \quad \mbox{for}\;\; t\geq 0.
\end{equation}

  {\em Proof}:\ \ Let the notation be as in Lemma 6.1. Since $e_k\geq 0$,
  by Lemmas 4.5 and 4.6 we see that $K(t)\geq 0$ and $\displaystyle{d\over dt}
  \big(e^{\kappa_0 t}K(t)\big)\leq 0$ for $t\geq 0$. By Lemma 6.2, it follows
  that the following estimate holds:
\begin{equation}
   |\Psi(t)|\leq |\tilde{\Psi}(t)|+K(0)\int^t_0e^{-\kappa_0(t-\tau)}
   |\tilde{\Psi}(\tau)|d\tau\;\;\;{\it for}\;\;t\geq 0.
\end{equation}
  Moreover, applying Lemma 4.5 to $e^{t\hat{\mathscr{L}}_k^+}e_k(r)$ and
  $e^{t\hat{\mathscr{L}}_k^+}\psi_0(r)$ we see that for any  nonnegative integer
  $k$ and any $\mu>\mu_k$ there exists corresponding constant $C=C_{k,\mu}>0$ such that
\begin{equation}
  \max_{0\leq r\leq R_s}|e^{t\hat{\mathscr{L}}_k^+}e_k(r)|\leq Ce^{\mu t}
  \quad \mbox{for}\;\; t\geq 0,
\end{equation}
\begin{equation}
  \max_{0\leq r\leq R_s}|e^{t\hat{\mathscr{L}}_k^+}\psi_0(r)|
  \leq Ce^{\mu t}\max_{0\leq r\leq R_s}|\psi_0(r)|
  \quad \mbox{for}\;\; t\geq 0.
\end{equation}
  The latter implies that
\begin{equation}
  |\tilde{\Psi}(t)|\leq Ce^{\mu t}\max_{0\leq r\leq R_s}|\psi_0(r)|
  \quad \mbox{for}\;\; t\geq 0.
\end{equation}
  Since $\mu_k\leq\mu_0^*$ for all nonnegative integer $k$, the above estimate holds for
  any $\mu>\mu_0^*$. Substituting (6.9) into (6.6) we easily see that for any nonnegative
  integer $k$ and any $\mu>\mu^*\equiv\max\{\mu_0^*,-\kappa_0\}$ there exists corresponding
  constant $C=C_{k,\mu}>0$ such that
\begin{equation}
  |\Psi(t)|\leq Ce^{\mu t}\max_{0\leq r\leq R_s}|\psi_0(r)|
  \quad \mbox{for}\;\; t\geq 0.
\end{equation}
  Now, noticing that
\begin{equation*}
  \hat{\mathscr{L}}_k\psi(r,t)=\hat{\mathscr{L}}_k^+\psi(r,t)-J(\psi(\cdot,t))e_k(r)
  =\hat{\mathscr{L}}_k^+\psi(r,t)-\Psi(t)e_k(r),
\end{equation*}
  by D'lHamul's formula we have
\begin{equation*}
  \psi(\cdot,t)=e^{t\hat{\mathscr{L}}_k^+}\psi_0
  -\int_0^t\Big(e^{(t-\tau)\hat{\mathscr{L}}_k^+}e_k\Big)\Psi(\tau)d\tau.
\end{equation*}
  From this relation and the estimates (6.7), (6.8) and (6.10), we immediately obtain
  (6.5). This proves Lemma 6.3. $\quad\Box$
\medskip

  {\bf Lemma 6.4}\ \ {\em Assume that the conditions in $(1.25)$ are satisfied and let $\kappa_0$ be as
  in Lemma 4.6. For every nonnegative integer $k$ there exists a corresponding constant
  $C=C_k>0$ such that for the solution of the problem $(6.3)$ the following estimate
  holds:}
\begin{equation}
  J(|\psi(\cdot,t)|)\leq CJ(|\psi_0|)(1+t)^2e^{-\kappa_0 t}
  \quad \mbox{for}\;\; t\geq 0.
\end{equation}

  {\em Proof}:\ \ Let the notation be as in Lemma 6.1. Applying Lemma 4.5 to
  $\tilde{\psi}(r,t)=e^{t\hat{\mathscr{L}}_k^+}\psi_0(r)$ we see that
$$
  \frac{d}{dt}J(|\tilde{\psi}(\cdot,t)|)\leq -\kappa_0 J(|\tilde{\psi}(\cdot,t)|)
  \quad \mbox{for}\;\; t\geq 0.
$$
  This implies that
$$
  |\tilde{\Psi}(t)|\leq J(|\tilde{\psi}(\cdot,t)|)\leq J(|\psi_0|)e^{-\kappa_0 t}
  \quad \mbox{for}\;\; t\geq 0.
$$
  From this estimate and Lemmas 6.1, 6.2 we get
\begin{equation}
  |\Psi(t)|\leq CJ(|\psi_0|)(1+t)e^{-\kappa_0 t}
  \quad \mbox{for}\;\; t\geq 0.
\end{equation}
  Now, we rewrite the equation $\displaystyle\partial_t\psi=\hat{\mathscr{L}}_k\psi$
  as follows:
$$
  \partial_t\psi(r,t)=\hat{\mathscr{L}}_k^+\psi(r,t)-\Psi(t)e_k(r).
$$
  Multiplying this equation with ${\rm sgn}\psi(r,t)$, we get
$$
  \partial_t|\psi(r,t)|\leq\hat{\mathscr{L}}_k^+|\psi(r,t)|+|\Psi(t)|e_k(r).
$$
  Using this relation and a similar argument as in the proof of Lemma 4.6 we get
\begin{equation}
  \frac{d}{dt}J(|\psi(\cdot,t)|)\leq -\kappa_0 J(|\psi(\cdot,t)|)+C|\Psi(t)|
  \quad \mbox{for}\;\; t\geq 0.
\end{equation}
  From (6.12) and (6.13) we easily see that (6.11) follows. $\quad\Box$
\medskip

  We are now ready to study the problem (5.1) for small $k$.
\medskip

  {\bf Lemma 6.5}\ \ {\em Assume that the conditions in $(1.25)$ are satisfied. There exists a constant
  $\mu^*<0$ such that for every nonnegative integer $k$ and any $\mu>\mu^*$ there
  exists corresponding constant $C=C_{k,\mu}>0$ such that for the solution of the
  problem $(5.1)$ the following estimate holds:}
\begin{equation}
  \max_{0\leq r\leq R_s}|\varphi(r,t)|\leq
  Ce^{\mu t}\max_{0\leq r\leq R_s}|\varphi_0(r)|
  \quad \mbox{for}\;\; t\geq 0.
\end{equation}

  {\em Proof}:\ \ Let $\varphi$ be the solution of the problem (5.1) and set
  $\psi(r,t)=r^{n+k-1}\varphi(r,t)/p_s'(r)$. By (6.2), $\psi$ is a solution of
  the problem (6.3) with initial data $\psi_0(r)=r^{n+k-1}\varphi_0(r)/p_s'(r)$.
  Hence, by Lemma 6.3 we see that for any $\mu>\max\{\mu_0^*,-\kappa_0\}$ there holds
\begin{equation*}
  \max_{0\leq r\leq R_s}r^{n+k-1}|\varphi(r,t)|/p_s'(r)\leq
  Ce^{\mu t}\max_{0\leq r\leq R_s}r^{n+k-1}|\varphi_0(r)|/p_s'(r)
  \quad \mbox{for}\;\; t\geq 0.
\end{equation*}
  This implies that for any $0<\delta<R_s$ there exists corresponding constant
  $C=C_{k,\mu,\delta}>0$ such that
\begin{equation}
  \max_{\delta\leq r\leq R_s}|\varphi(r,t)|\leq
  Ce^{\mu t}\max_{0\leq r\leq R_s}|\varphi_0(r)|
  \quad \mbox{for}\;\; t\geq 0.
\end{equation}
  Leaving $\delta$ to be specified later, we take a nonnegative cut-off function
  $\chi\in C[0,R_s]$ such that
$$
  \chi(r)\leq 1 \;\; \mbox{for}\;\; 0\leq r\leq R_s, \quad
  \chi(r)=1 \;\; \mbox{for}\;\;  0\leq r\leq\delta, \quad \mbox{and} \quad
  \chi(r)=0 \;\; \mbox{for}\;\; 2\delta\leq r\leq R_s,
$$
  and split the operator $\mathscr{K}_k$ introduced in (5.4) into the sum of
  two operators $\mathscr{K}_k'$ and $\mathscr{K}_k''$, where
\begin{eqnarray*}
  \mathscr{K}_k'(\phi)&=&r^{k-1}p_s'(r)\Big[{\theta_k}
  \int_r^{\max\{r,\delta\}}\!\!\!\rho^{-k+1}g_p^*(\rho)\phi(\rho)d\rho
  -\frac{(1\!-\!\theta_k)\chi(r)}{r^{n+2(k-1)}}
  \int_0^{\min\{r,\delta\}}\!\!\!\rho^{n+k-1}g_p^*(\rho)\phi(\rho)d\rho\Big],
\\
  \mathscr{K}_k''(\phi)&=&r^{k-1}p_s'(r)\Big[{\theta_k}
  \int_{\max\{r,\delta\}}^{R_s}\!\rho^{-k+1}g_p^*(\rho)\phi(\rho)d\rho
  +\frac{1-\theta_k}{R_s^{n+2(k-1)}}
  \int_0^{R_s}\!\!\!\rho^{n+k-1}g_p^*(\rho)\phi(\rho)d\rho
\nonumber\\ [0.3cm]
  &&\displaystyle-\frac{(1\!-\!\theta_k)[1-\chi(r)]}{r^{n+2(k-1)}}
  \int_0^{\min\{r,\delta\}}\!\!\!\rho^{n+k-1}g_p^*(\rho)\phi(\rho)d\rho
  -\frac{1\!-\!\theta_k}{r^{n+2(k-1)}}
  \int_{\min\{r,\delta\}}^r\!\!\!\rho^{n+k-1}g_p^*(\rho)\phi(\rho)d\rho\Big].
\end{eqnarray*}
  Next we let $f(r,t)=\mathscr{K}_k''(\varphi(\cdot,t))(r)$. Since
  $\tilde{\mathscr{L}}_k=\mathscr{L}_0+\mathscr{K}_k$, from (5.1) we see that
  $\varphi$ is the solution of the following problem:
\begin{equation}
\left\{
\begin{array}{l}
  \partial_t\varphi=\mathscr{L}_0(\varphi)+\mathscr{K}_k'(\varphi)
  +f(r,t) \quad \mbox{for}\;\; 0\leq r\leq 1,\;\; t>0,\\
  \varphi|_{t=0}=\varphi_0\;\; \mbox{for}\;\; 0\leq r\leq R_s.
\end{array}
\right.
\end{equation}
  By using (6.11) and (6.15) we easily see that for any $\mu>\max\{\mu_0^*,-\kappa_0\}$
  there exists corresponding constant $C=C(k,\mu,\delta)>0$ such that
\begin{equation}
  \max_{0\leq r\leq R_s}|f(r,t)|\leq Ce^{\mu t}\max_{0\leq r\leq R_s}|\varphi_0(r)|
  \quad \mbox{for}\;\; t\geq 0.
\end{equation}
  By Corollary 4.3 we see that for any $\mu>\nu_0$ there holds
\begin{equation}
  \max_{0\leq r\leq R_s}|e^{t\mathscr{L}_0}\phi(r)|\leq C_{\mu}
  e^{\mu t}\max_{0\leq r\leq R_s}|\phi(r)|
  \quad \mbox{for}\;\; t\geq 0.
\end{equation}
  Besides, it is easy to check that there exists a positive function $\varepsilon(\delta)$
  of $\delta$ which converges to zero as $\delta\to 0^+$, such that
\begin{equation}
  \max_{0\leq r\leq R_s}|\mathscr{K}_k'(\phi)|\leq\varepsilon(\delta)
  \max_{0\leq r\leq R_s}|\phi(r)|
  \quad \mbox{for}\;\; \phi\in C[0,R_s].
\end{equation}
  For instance, for the cases $k\geq 2$ we have
\begin{eqnarray*}
  &&\max_{0\leq r\leq R_s}\Big|r^{k-1}p_s'(r)\int_r^{\max\{r,\delta\}}\!\!\!
  \rho^{-k+1}g_p^*(\rho)\phi(\rho)d\rho\Big|
\\
  &\leq& C\max_{0\leq r\leq\delta}\Big(rp_s'(r)\int_r^{\delta}
  \rho^{-1}d\rho\Big)\cdot\max_{0\leq r\leq R_s}|\phi(r)|
\\
  &=& C\varepsilon_1(\delta)\max_{0\leq r\leq R_s}|\phi(r)|,
\end{eqnarray*}
  where $\varepsilon_1(\delta)=\displaystyle\max_{0\leq r\leq\delta}\Big(rp_s'(r)\ln
  \frac{\delta}{r}\Big)\to 0$ as $\delta\to 0^+$. For the cases $k=0,1$ we can use
  Lemma 3.2 to get a similar inequality. By a standard perturbation theorem for
  $C_0$-semigroups, from (6.18) and (6.19) we have
\begin{equation}
  \max_{0\leq r\leq R_s}|e^{t(\mathscr{L}_0+\mathscr{K}_k')}\phi(r)|\leq C_{\mu}'
  e^{[\mu+C_{\mu}\varepsilon(\delta)]t}\max_{0\leq r\leq R_s}|\phi(r)|
  \quad \mbox{for}\;\; t\geq 0.
\end{equation}
  Now, from (6.16) we have
\begin{equation}
  \varphi(\cdot,t)=e^{t(\mathscr{L}_0+\mathscr{K}_k')}\varphi_0+
  \int_0^te^{(t-\tau)(\mathscr{L}_0+\mathscr{K}_k')}f(\cdot,\tau)d\tau
  \quad \mbox{for}\;\; t\geq 0.
\end{equation}
  From (6.17), (6.20) and (6.21) we can easily deduce that for any given $\mu>\mu^*
  \equiv\max\{\mu_0^*,-\kappa_0\}$, by first choosing $\mu'>\mu^*$ such that $\mu>\mu'$
  and next choosing $\delta$ sufficiently small so that $\mu'+C_{\mu}\varepsilon(\delta)
  <\mu$, (6.14) follows. This completes the proof. $\quad\Box$
\medskip

  {\bf Lemma 6.6}\ \ {\em  Assume that the conditions in $(1.25)$ are satisfied. There exists a constant
  $\mu_0^*<0$ such that for every nonnegative integer $k$, any $\alpha\geq 1$ and $\mu
  >\mu_{\alpha}^*$ there exist corresponding positive constant $C=C_{k,\alpha,\mu}$
  such that for any initial data $\varphi_0$ such that $\varphi_0\in L^{\alpha}([0,R_s];
  r^{n-1}dr)$, the solution of the initial value problem $(5.1)$ satisfies the following
  estimate:}
\begin{equation}
  \Big(\int_0^{R_s}|\varphi(r,t)|^{\alpha}r^{n-1}dr\Big)^{\frac{1}{\alpha}}
  \leq C\Big(\int_0^{R_s}|\varphi_0(r)|^{\alpha}r^{n-1}dr\Big)^{\frac{1}{\alpha}}e^{\mu t}
  \quad \mbox{for}\;\; t\geq 0.
\end{equation}

  {\em Proof}:\ \ We split the proof into three steps.

  {\em Step 1}:\ \ We first prove that
\begin{equation}
  \int_0^{R_s}|\varphi(r,t)|r^{n+k-1}dr\leq Ce^{\mu t}\int_0^{R_s}|\varphi_0(r)|r^{n+k-1}dr
  \quad \mbox{for}\;\; t\geq 0.
\end{equation}
  Indeed, by the relation (6.2) it follows that if $\varphi$ is a solution
  of the problem (5.1) then $\psi(r,t)=r^{n+k-1}\varphi(r,t)/p_s'(r)$ is a solution
  of the problem (6.3) with initial data $\psi_0(r)=r^{n+k-1}\varphi_0(r)/p_s'(r)$.
  Due to this fact, the above estimate is an immediate consequence of Lemma 6.4.

  {\em Step 2}:\ \ We next prove that
\begin{equation}
  \int_0^{R_s}|\varphi(r,t)|r^{n-1}dr\leq Ce^{\mu t}\int_0^{R_s}|\varphi_0(r)|r^{n-1}dr
  \quad \mbox{for}\;\; t\geq 0.
\end{equation}
  Let $\mathscr{K}_k'$, $\mathscr{K}_k''$ and $f$ be as in the proof of Lemma 6.5.
  Using (6.23) we easily see that
\begin{equation}
  \int_0^{R_s}|f(r,t)|r^{n-1}dr\leq C_{\delta,\mu}e^{\mu t}\int_0^{R_s}|\varphi_0(r)|r^{n-1}dr
  \quad \mbox{for}\;\; t\geq 0.
\end{equation}
  By Lemma 4.3 we see that with $\nu^*=\displaystyle\max_{0\leq r\leq R_s}
  [g^*(r)+f_p^*(r)]<0$ there holds
\begin{equation}
  \int_0^{R_s}|e^{t\mathscr{L}_0}\phi(r)|r^{n-1}dr\leq
  e^{\nu^* t}\int_0^{R_s}|\varphi_0(r)|r^{n-1}dr
  \quad \mbox{for}\;\; t\geq 0.
\end{equation}
  Besides, we have
\begin{eqnarray}
  \int_0^{R_s}|\mathscr{K}_k'(\phi)|r^{n-1}dr
  &\leq&\theta_k\int_0^{\delta}\int_r^{\delta}r^{n+k-2}p_s'(r)
  \rho^{-k+1}g_p^*(\rho)|\phi(\rho)|d\rho dr
\nonumber\\
  &&+(1\!-\!\theta_k)\int_0^{2\delta}\int_0^r r^{-k}p_s'(r)
  \rho^{n+k-1}g_p^*(\rho)|\phi(\rho)|d\rho dr
\nonumber\\
  &\leq&C\int_0^{\delta}\Big(\int_0^{\rho}r^{n+k-2}p_s'(r)dr\Big)
  \rho^{-k+1}|\phi(\rho)|d\rho
\nonumber\\
  &&+C\int_0^{2\delta}\Big(\int_{\rho}^{2\delta}r^{-k}p_s'(r)dr\Big)
  \rho^{n+k-1}|\phi(\rho)|d\rho
\nonumber\\
  &\leq&C[p_s(\delta)-p_s(0)]\int_0^{\delta}\rho^{n-1}|\phi(\rho)|d\rho
  +C[p_s(2\delta)-p_s(0)]\int_0^{2\delta}\rho^{n-1}|\phi(\rho)|d\rho
\nonumber\\
  &\leq&\varepsilon(\delta)\int_0^{R_s}\rho^{n-1}|\phi(\rho)|d\rho,
\end{eqnarray}
  where $\varepsilon(\delta)=2C[p_s(2\delta)-p_s(0)]$. By using (6.21), (6.25),
  (6.26), (6.27) and a similar argument as in the proof of Lemma 6.5, we
  obtain (6.24).

  {\em Step 3}:\ \   Interpolating the inequalities (6.21) and (6.24), we see
  that (6.22) follows. $\quad\Box$

\section{Decay estimates for the system (2.22)--(2.23)}
\setcounter{equation}{0}

\hskip 2em
  {\bf Lemma 7.1}\ \ {\em Assume that the conditions in $(1.25)$ are satisfied. There exists constants
  $\gamma^*>0$ and $\lambda>0$ such that for any $\gamma>\gamma^*$, any
  integer $k\geq 2$ and any $\alpha\geq 1$, the solution of the system of equations
  $(2.22)$--$(2.23)$ satisfies the following estimates:
\begin{equation}
  \max_{0\leq r\leq R_s}|\varphi_k(r,t)|+|\eta_k(t)|\leq
  Ce^{-\lambda t}[\max_{0\leq r\leq R_s}|\varphi_{k0}(r)|+|\eta_{k0}|]
  \quad \mbox{for}\;\; t\geq 0,
\end{equation}
\begin{equation}
  \Big(\int_0^{R_s}|\varphi_k(r,t)|^{\alpha}r^{n-1}dr\Big)^{\frac{1}{\alpha}}
  +|\eta_k(t)|\leq C_{\alpha}e^{-\lambda t}\Big[
  \Big(\int_0^{R_s}|\varphi_{k0}(r)|^{\alpha}r^{n-1}dr\Big)^{\frac{1}{\alpha}}
  +|\eta_{k0}|\Big] \quad \mbox{for}\;\; t\geq 0,
\end{equation}
  where $\varphi_{k0}$ and $\eta_{k0}$ are the initial data of $\varphi_k$ and
  $\eta_k$, respectively, and $C$, $C_{\alpha}$ represent positive constants
  independent of $k$.}
\medskip

  {\em Proof}:\ \ We only need to prove that the solution of the system of equations
  $(2.24)$--$(2.25)$ satisfies the following estimates:
\begin{equation}
  \max_{0\leq r\leq R_s}|\tilde{\varphi}_k(r,t)|+|\eta_k(t)|\leq
  Ce^{-\lambda t}[\max_{0\leq r\leq R_s}|\tilde{\varphi}_{k0}(r)|+|\eta_{k0}|]
  \quad \mbox{for}\;\; t\geq 0,
\end{equation}
\begin{equation}
  \Big(\int_0^{R_s}|\tilde{\varphi}_k(r,t)|^{\alpha}r^{n-1}dr\Big)^{\frac{1}{\alpha}}
  +|\eta_k(t)|\leq C_{\alpha}e^{-\lambda t}\Big[
  \Big(\int_0^{R_s}|\tilde{\varphi}_{k0}(r)|^{\alpha}r^{n-1}dr\Big)^{\frac{1}{\alpha}}
  +|\eta_{k0}|\Big] \quad \mbox{for}\;\; t\geq 0,
\end{equation}
  where $\tilde{\varphi}_{10}$ and $\eta_{10}$ are the initial data of
  $\tilde{\varphi}_1$ and $\eta_1$, respectively. Indeed, since the solutions of the
  systems $(2.22)$--$(2.23)$ and $(2.24)$--$(2.25)$ satisfy the relations
\begin{equation}
  \widetilde{\varphi}_k(r,t)=\varphi_k(r,t)+\Big(\frac{r}{R_s}\Big)^{k-1}p_s'(r)
  \eta_k(t), \quad k=0,1,2,\cdots,
\end{equation}
  we see that (7.1) and (7.2) are immediate consequences of (7.3) and (7.4).

  Fix two constants $\lambda$ and $\mu$ such that $\lambda>0$
  and $\mu^*<\mu<-\lambda$, where $\mu^*$ is as
  in Lemma 6.5. We re-denote the constant appearing on the right-hand side of
  (6.14) as $C_0$. Let $\tilde{\varphi}_{k0}\in C[0,R_s]$ and $\eta_{k0}\in
  \mathbb{R}$ be given. Given $\varphi\in C([0,R_s]\times[0,\infty))$ satisfying the
  condition
\begin{equation}
  \max_{0\leq r\leq R_s}|\varphi(r,t)|\leq 2C_0e^{-\lambda t}
  [\max_{0\leq r\leq R_s}|\tilde{\varphi}_{k0}(r)|+|\eta_{k0}|]
  \quad \mbox{for}\;\; t\geq 0,
\end{equation}
  we consider the following initial value problems:
\begin{eqnarray}
  &&\frac{d\eta}{dt}=\tilde{\alpha}_k(\gamma)\eta
  +J_k(\varphi)\quad \mbox{for}\;\; t\geq 0,\quad
  \mbox{and}\;\; \eta|_{t=0}=\eta_{k0},
\\
  &&\frac{\partial\tilde{\varphi}}{\partial t}
  =\tilde{\mathscr{L}}_k(\widetilde{\varphi})
  +c_k(r)\eta\quad \mbox{for}\;\; 0<r<R_s,\;\; t\geq 0,\quad
  \mbox{and}\;\; \tilde{\varphi}|_{t=0}=\tilde{\varphi}_{k0}.
\end{eqnarray}
  The solution of (7.7) is given by
\begin{equation*}
  \eta(t)=e^{\tilde{\alpha}_k(\gamma)t}\eta_{k0}+
  \int_0^te^{\tilde{\alpha}_k(\gamma)(t-\tau)}J_k(\varphi(\cdot,\tau))d\tau
  \quad \mbox{for}\;\; t\geq 0.
\end{equation*}
  It is clear that there exists a constant $c>0$ independent of $k$ and $\gamma$
  such that for $\gamma$ sufficiently large,
\begin{equation}
  \tilde{\alpha}_k(\gamma)\leq -ck^3\gamma \quad \mbox{for}\;\; k\geq 2.
\end{equation}
  Using this fact and (7.6) we easily see that for $\gamma$ sufficiently large and
  $k\geq 2$,
\begin{equation}
  |\eta(t)|\leq e^{-ck^3\gamma t}|\eta_{k0}|+
  \frac{CC_0e^{-\lambda t}}{ck^3\gamma-\lambda}
  [\max_{0\leq r\leq R_s}|\tilde{\varphi}_{k0}(r)|+|\eta_{k0}|]
  \quad \mbox{for}\;\; t\geq 0.
\end{equation}
  Having solved the initial value problem (7.7), we substitute its solution into
  (7.8) and next solve that initial value problem. The solution is given by
\begin{equation}
  \tilde{\varphi}(\cdot,t)=e^{t\tilde{\mathscr{L}}_k}\tilde{\varphi}_{k0}+
  \int_0^t[e^{(t-\tau)\tilde{\mathscr{L}}_k}c_k]\eta(\tau)d\tau
  \quad \mbox{for}\;\; t\geq 0.
\end{equation}
  It is clear that there exists constant $C>0$ independent of $k$ such that
  for any $k\geq 2$,
\begin{equation}
  \max_{0\leq r\leq R_s}|c_k(r)|\leq Ck.
\end{equation}
  By using Lemmas 5.1, 6.5 and the estimates (7.10), (7.12) we easily deduce from
  (7.11) to get
\begin{eqnarray*}
  \max_{0\leq r\leq R_s}|\tilde{\varphi}(r,t)|
  &\leq & C_0e^{\mu t}\max_{0\leq r\leq R_s}|\tilde{\varphi}_{k0}(r)|
  +\frac{CC_0ke^{-\lambda t}}{ck^3\gamma+\mu}|\eta_{k0}|
\nonumber\\
  &&+\frac{CC_0ke^{-\lambda t}}{|\lambda+\mu|(ck^3\gamma-\lambda)}
  [\max_{0\leq r\leq R_s}|\tilde{\varphi}_{k0}(r)|+|\eta_{k0}|].
\end{eqnarray*}
  From this estimate we easily see that by choosing $\gamma^*>0$ sufficiently large,
  we have that for any $\gamma>\gamma^*$ and any $k\geq 2$,
\begin{equation*}
  \max_{0\leq r\leq R_s}|\tilde{\varphi}(r,t)|\leq  2C_0e^{-\lambda t}
  [\max_{0\leq r\leq R_s}|\tilde{\varphi}_{k0}(r)|+|\eta_{k0}|]
  \quad \mbox{for}\;\; t\geq 0,
\end{equation*}
  i.e., $\tilde{\varphi}$ satisfies the condition (7.6). Furthermore, it can also
  be easily seen that if we choose $\gamma^*>0$ so large that for any $\gamma>\gamma^*$
  and any $k\geq 2$,
$$
  \displaystyle\frac{Ck}{|\lambda+\mu|(ck^3\gamma-\lambda)}<1,
$$
  then the mapping $\varphi\mapsto\tilde{\varphi}$ is a contraction. Hence, by using
  the standard contraction mapping argument we see that the system of equations
  (2.24)--(2.25) subject to the initial conditions in (7.7) and (7.8) has a unique
  solution satisfying (7.6) (with $\varphi$ replaced by $\tilde{\varphi}$) and (7.10).
  This proves the estimate (7.3). To prove the estimate (7.4), we replace the
  condition (7.6) with the following one:
\begin{equation*}
  \Big(\int_0^{R_s}|\varphi_k(r,t)|^{\alpha}r^{n-1}dr\Big)^{\frac{1}{\alpha}}
  \leq 2C_0e^{-\lambda t}\Big[
  \Big(\int_0^{R_s}|\tilde{\varphi}_{k0}(r)|^{\alpha}r^{n-1}dr\Big)^{\frac{1}{\alpha}}
  +|\eta_{k0}|\Big] \quad \mbox{for}\;\; t\geq 0,
\end{equation*}
  and use a similar argument as above but instead of using Lemmas 5.1 and 6.5 we now
  use Lemmas 5.2 and 6.6. We omit the details. This completes the proof of Lemma 7.1.
  $\quad\Box$
\medskip

  The above results do not work for the cases $k=0,1$. We first note that in these
  cases $\widetilde{\alpha}_k$ are independent of $\gamma$. For these special cases,
  we have the following results:

  {\bf Lemma 7.2}\ \ {\em In the case $k=1$, there exist constants $\lambda>0$ and
  $C>0$ such that the solution of the system of equations $(2.24)$--$(2.25)$
  satisfies the following estimates:
\begin{equation}
  \max_{0\leq r\leq R_s}|\varphi_1(r,t)-\varphi_{\infty}(r)|\leq
  Ce^{-\lambda t}\Big(\max_{0\leq r\leq R_s}|\varphi_{10}(r)|
  +|\eta_{10}|\Big) \quad \mbox{for}\;\; t\geq 0,
\end{equation}
\begin{eqnarray}
  \Big(\int_0^{R_s}|\varphi_1(r,t)-\varphi_{\infty}(r)|^{\alpha}r^{n-1}dr
  \Big)^{\frac{1}{\alpha}}\leq && C_{\alpha}e^{-\lambda t}
  \Big(\int_0^{R_s}|\tilde{\varphi}_{10}(r)|^{\alpha}r^{n-1}dr
  +|\eta_{10}|\Big)^{\frac{1}{\alpha}}
\nonumber\\
   && \qquad \mbox{for}\;\; \alpha\geq 1,\;\; t\geq 0,
\end{eqnarray}
\begin{equation}
  |\eta_1(t)-\eta_{\infty}|\leq Ce^{-\lambda t}
  \Big(\max_{0\leq r\leq R_s}|\varphi_{10}(r)|+|\eta_{10}|\Big)
  \quad \mbox{for}\;\; t\geq 0,
\end{equation}
  where $\varphi_{\infty}(r)=-p_s'(r)\eta_{\infty}$, and $\eta_{\infty}$ is a real
  constant uniquely determined by the initial data $\varphi_{10}$ and $\eta_{10}$.}
\medskip

  {\em Proof}:\ \ Because of the relations in (7.5) and Lemma 3.2, we see that the
  above estimates follow if we prove that the solution of the system $(2.24)$--$(2.25)$
  satisfies the following estimates:
\begin{equation}
  \max_{0\leq r\leq R_s}|\tilde{\varphi}_1(r,t)|\leq
  Ce^{-\lambda t}\max_{0\leq r\leq R_s}|\tilde{\varphi}_{10}(r)|
  \quad \mbox{for}\;\; t\geq 0,
\end{equation}
\begin{equation}
  \Big(\int_0^{R_s}|\tilde{\varphi}_1(r,t)|^{\alpha}r^{n-1}dr\Big)^{\frac{1}{\alpha}}
  \leq C_{\alpha}e^{-\lambda t}
  \Big(\int_0^{R_s}|\tilde{\varphi}_{10}(r)|^{\alpha}r^{n-1}dr\Big)^{\frac{1}{\alpha}}
   \quad \mbox{for}\;\; \alpha\geq 1,\;\; t\geq 0,
\end{equation}
\begin{equation}
  \Big|\eta_1(t)-\eta_{10}-\int_0^{\infty}J_1(\tilde{\varphi}_1(\cdot,t))dt\Big|
  \leq Ce^{-\lambda t}\max_{0\leq r\leq R_s}|\tilde{\varphi}_{10}(r)|
  \quad \mbox{for}\;\; t\geq 0.
\end{equation}
  These estimates are immediate consequences of the relations in (2.29) and
  Lemmas 6.4--6.6 applied to the case $k=1$. $\quad\Box$
\medskip

  {\bf Lemma 7.3}\ \ {\em In the case $k=0$, there exist constants $\lambda>0$ and
  $C>0$ such that the solution of the system of equations $(2.22)$--$(2.23)$
  satisfies the following estimate:
\begin{equation}
  \max_{0\leq r\leq R_s}|\varphi_0(r,t)|+|\eta_0(t)|\leq
  Ce^{-\lambda t}[\max_{0\leq r\leq R_s}|\varphi_{00}(r)|+|\eta_{00}|]
  \quad \mbox{for}\;\; t\geq 0,
\end{equation}
\begin{equation}
  \Big(\int_0^{R_s}|\varphi_0(r,t)|^{\alpha}r^{n-1}dr\Big)^{\frac{1}{\alpha}}
  +|\eta_0(t)|\leq C_{\alpha}e^{-\lambda t}\Big[
  \Big(\int_0^{R_s}|\varphi_{00}(r)|^{\alpha}r^{n-1}dr\Big)^{\frac{1}{\alpha}}
  +|\eta_{00}|\Big] \quad \mbox{for}\;\; t\geq 0,
\end{equation}
  where $\varphi_{00}$ and $\eta_{00}$ are the initial data of $\varphi_0$ and $\eta_0$,
  respectively.}
\medskip

  {\em Proof}:\ \ In the case $n=3$, the estimate (7.19) follows from Lemma 6.2 of
  \cite{Cui2}, which is an improvement of Theorem 5.1 of \cite{ChenCuiF}. This is
  because the system (2.22)--(2.23) is the linearization of the radial version of the
  system (1.12)--(1.18), so that it is equivalent to the system (5.1)--(5.2) in
  \cite{ChenCuiF}. The estimate (7.20) follows from this fact and a similar argument
  as in the proof of Lemma 6.6 and (7.2) in Lemma 7.1. In the general dimension case,
  the arguments are similar. In what follows we give a sketch of them. Since
  $b_0(r,\gamma)$ and $\alpha_0(\gamma)$ are actually independent of $\gamma$, in the
  sequel we briefly write them as $b_0(r)$ and $\alpha_0$, respectively.

  Firstly, by using some similar arguments as in the proofs of Lemmas 7.1--7.3 of
  \cite{ChenCuiF} (see also the proof of Lemma 4.4 of this work) we can prove that the
  equation
\begin{equation}
  \mathscr{L}_0\phi^*(r)-[\alpha_0+J_0(\phi^*)]\phi^*(r)+b_0(r)=0
  \quad \mbox{for}\;\; 0<r<R_s
\end{equation}
  has a unique solution $\phi^*(r)=\psi^*(r)-\displaystyle\Big(\frac{R_s}{r}\Big)^{n-1}
  p_s'(r)$ with $\alpha_0+J_0(\phi^*)<0$, where $\psi^*\in C^1(0,R_s]$, $\psi^*(r)>0$
  for $0<r<R_s$, and $J_0(\psi^*)<\infty$. More precisely, the above equation can be
  rewritten into the following equivalent system of equations:
\begin{equation*}
\left\{
\begin{array}{l}
  [\lambda-\mathscr{L}_0^+]\psi^*(r)=h_0(r) \quad \mbox{for}\;\; 0<r<R_s
\\ [0.1cm]
  \lambda=\displaystyle\alpha_0+J_0(\psi^*)-J_0\Big[\Big(\frac{R_s}{r}\Big)^{n-1}p_s'(r)\Big]
\end{array}
\right.
\end{equation*}
  (with $\psi^*$ and $\lambda$ being the unknowns), where
\begin{equation*}
  \mathscr{L}_0^+\varphi(r)=-v_s(r)\varphi'(r)+f_p^*(r)\varphi(r)
  +\frac{p_s'(r)}{r^{n-1}}\int_r^{R_s}\rho^{n-1}f_p^*(\rho)\varphi(\rho)d\rho,
\end{equation*}
\begin{equation*}
  h_0(r)=b_0(r)-\mathscr{L}_0^+\Big[\Big(\frac{R_s}{r}\Big)^{n-1}p_s'(r)\Big]
  +\alpha_0\Big(\frac{R_s}{r}\Big)^{n-1}p_s'(r).
\end{equation*}
  By using the assertions (4) and (5) of Lemma 3.3 and the equations (2.3), (2.4) and
  (2.5) as well as their equivalent integral forms, we can prove that $h_0(r)>0$ for
  $0<r<R_s$. Hence, by using some similar arguments as in the proofs of Lemmas 7.1--7.3
  of \cite{ChenCuiF}, we see that the above system has a unique solution $(\psi^*,
  \lambda)$ with the properties that
\begin{equation*}
  nv_s'(0)+f_p^*(0)<\lambda<0, \quad  \lambda\geq f_p^*(R_s),
\end{equation*}
  and
\begin{equation*}
  c_1r^{-(n-1)}(R_s-r)p_s'(r)\leq \psi^*(r)\leq c_2r^{-\theta}\;\; \mbox{for}\;\; 0<r<R_s,
\end{equation*}
  where $c_1$, $c_2$ and $\theta$ are positive constants, $0<\theta<n$, so that
  $J_0(\psi^*)<\infty$. This proves the above statement. We now make a transformation
  of unknown variables $(\varphi_0,\eta_0)\mapsto (\psi_0,\eta_0)$ such that
\begin{equation}
  \psi_0(r,t)=\varphi_0(r,t)-\phi^*(r)\eta(t).
\end{equation}
  Then the system (2.22)--(2.23) with $k=0$ is transformed into the following
  equivalent system:
\begin{eqnarray}
  \frac{\partial\psi_0}{\partial t}&=&\mathscr{L}_0(\psi_0)-\phi^*(r)J_0(\psi_0),
\\
  \frac{d\eta_0}{dt}&=&\lambda\eta_0+J_0(\psi_0).
\end{eqnarray}
  Note that
\begin{equation}
  \mathscr{L}_0(\psi_0)-\phi^*(r)J_0(\psi_0)
  =\mathscr{L}_0^+(\psi_0)-\psi^*(r)J_0(\psi_0).
\end{equation}
  Using this fact and some similar arguments as in the proofs of Lemmas 6.5 and
  6.6, we obtain from (7.23) the following estimates:
\begin{equation}
  \max_{0\leq r\leq R_s}|\psi_0(r,t)|\leq
  Ce^{-\lambda t}\max_{0\leq r\leq R_s}|\psi_0(r,0)|
  \quad \mbox{for}\;\; t\geq 0,
\end{equation}
\begin{equation}
  \Big(\int_0^{R_s}|\psi_0(r,t)|^{\alpha}r^{n-1}dr\Big)^{\frac{1}{\alpha}}
  \leq C_{\alpha}e^{-\lambda t}
  \Big(\int_0^{R_s}|\psi_0(r,0)|^{\alpha}r^{n-1}dr\Big)^{\frac{1}{\alpha}}
  \quad \mbox{for}\;\; t\geq 0.
\end{equation}
  By (7.24), this further implies that
\begin{equation}
  |\eta_0(t)|\leq Ce^{-\lambda t}\Big[
  \int_0^{R_s}|\psi_0(r,0)|r^{n-1}dr+|\eta_{00}|\Big]
  \quad \mbox{for}\;\; t\geq 0.
\end{equation}
  The estimates (7.19) and (7.20) now follow from (7.22) and (7.26)--(7.28).
  $\quad\Box$
\medskip

\section{Decay estimates for the system (2.19)}
\setcounter{equation}{0}

\hskip 2em
  For every nonnegative integer $k$, we let $Y_{kl}(\omega)$ ($l=1,2,\cdots,d_k$) be
  a normalized orthogonal basis of the linear space of spherical harmonics of degree
  $k$ (cf. \cite{SW}), i.e.
$$
  \Delta_\omega Y_{kl}(\omega)=-\lambda_k Y_{kl}(\omega), \quad
  \lambda_k=(n+k-2)k, 
$$
$$
  \int_{\mathbb{S}^{n-1}}Y_{kl}(\omega)Y_{kl'}(\omega)d\omega=0\;\; (l\neq l'),
  \qquad
  \int_{\mathbb{S}^{n-1}}Y_{kl}^2(\omega)d\omega=1.
$$
  Here $d_k$ is the dimension of the linear space of spherical harmonics of degree
  $k$, i.e.
\begin{equation}
  d_0=1, \quad d_1=n \quad \mbox{and} \quad
  d_k={n\!+\!k\!-\!1\choose k}-{n\!+\!k\!-\!3\choose k\!-\!2}
  \quad\hbox{for}\quad k\geq 2,
\end{equation}
  and $d\omega$ is the induced element on the unit sphere $\mathbb{S}^{n-1}$ of the
  Lebesque measure $dx$ in $\mathbb{R}^n$. Note that in particular,
$$
  Y_{01}(\omega)=\frac{1}{\sqrt{\sigma_n}} \quad \mbox{and} \quad
  Y_{1l}(\omega)=\frac{\sqrt{n}\omega_l}{\sqrt{\sigma_n}}, \quad l=1,2,\cdots,n,
$$
  where $\sigma_n$ denotes the surface area of the unit sphere $\mathbb{S}^{n-1}$,
  i.e. $\sigma_n=\displaystyle\frac{2\pi^{n/2}}{\Gamma(n/2)}$, and $\omega_l$
  denotes the $l$-th component of $\omega$.

  For any $1\leq\alpha<\infty$ and $1\leq\beta<\infty$, we denote by
  $X_{\alpha\beta}$ the space of all measurable functions $u(x)$ in the ball
  $\mathbb{B}(0,R_s)\subseteq\mathbb{R}^n$ satisfying the following conditions:
\begin{equation}
  u(x)=\sum_{k=0}^{\infty}\sum_{l=1}^{d_k}u_{kl}(r)Y_{kl}(\omega)\;\;
  \mbox{in}\;\; S'(\mathbb{B}(0,R_s)), \quad
  r=|x|,\;\; \omega=\frac{x}{|x|},
\end{equation}
$$
  \|u\|_{X_{\alpha\beta}}=\Big[\sum_{k=0}^{\infty}\sum_{l=1}^{d_k}
  \Big(\int_0^{R_s}|u_{kl}(r)|^{\alpha}r^{n-1}dr\Big)^{\frac{\beta}{\alpha}}
  \Big]^{\frac{1}{\beta}}<\infty.
$$
  We also let $X_{\infty\beta}$, $X_{\alpha\infty}$ and $X_{\infty\infty}$ be
  respectively the following spaces:
$$
  \|u\|_{X_{\infty\beta}}=\Big[\sum_{k=0}^{\infty}\sum_{l=1}^{d_k}
  \Big(\sup_{0\leq r\leq R_s}|u_{kl}(r)|\Big)^{\beta}\Big]^{\frac{1}{\beta}}<\infty,
$$
$$
  \|u\|_{X_{\alpha\infty}}=\sup_{k,l}
  \Big(\int_0^{R_s}|u_{kl}(r)|^{\alpha}r^{n-1}dr\Big)^{\frac{1}{\alpha}}<\infty,
$$
$$
  \|u\|_{X_{\infty\infty}}=\sup_{k,l}\sup_{0\leq r\leq R_s}|u_{kl}(r)|<\infty.
$$
  It is clear that for any $1\leq\alpha\leq\infty$ and $1\leq\beta\leq\infty$,
  $X_{\alpha\beta}$ is a Banach space. For $1\leq\beta\leq\infty$, we denote by
  $\dot{X}_{\infty\beta}$ the closure of $C(\overline{\mathbb{B}(0,R_s)})$ in
  $X_{\infty\beta}$. If $u\in\dot{X}_{\infty\beta}$ then from the relation
$$
  u_{kl}(r)=\int_{\mathbb{S}^{n-1}}u(r\omega)Y_{kl}(\omega)d\omega,
  \quad k=0,1,2,\cdots,\;\; l=1,2,\cdots,d_k,
$$
  we see that in the expansion (8.2), all coefficients $u_{kl}(r)$ are
  continuous functions in $[0,R_s]$. Note that $X_{22}=L^2(\mathbb{B}(0,R_s))$.

  Next, for any $1\leq\beta<\infty$, we denote by $Y_{\beta}$ the space of all
  measurable functions $\varphi(\omega)$ on the sphere $\mathbb{S}^{n-1}$
  satisfying the following conditions:
\begin{equation}
  \varphi(\omega)=\sum_{k=0}^{\infty}\sum_{l=1}^{d_k}a_{kl}Y_{kl}(\omega)\;\;
  \mbox{in}\;\; D'(\mathbb{S}^{n-1}),
\end{equation}
$$
  \|\varphi\|_{Y_{\beta}}=\Big(\sum_{k=0}^{\infty}\sum_{l=1}^{d_k}|a_{kl}|^{\beta}
  \Big)^{\frac{1}{\beta}}<\infty.
$$
  We also denote by $Y_{\infty}$ the following space:
$$
  \|\varphi\|_{Y_{\infty}}=\sup_{k,l}|a_{kl}|<\infty.
$$
  It is clear that for any $1\leq\beta\leq\infty$, $Y_{\beta}$ is a Banach space.
  Note that $Y_2=L^2(\mathbb{S}^{n-1})$.

  We can now give the precise statement of the main result of this paper:
\medskip

  {\bf Theorem 8.1}\ \ {\em Assume that the conditions in $(1.25)$ are satisfied.
  There exists constants $\gamma^*>0$ and $\lambda>0$ such that for any $\gamma>
  \gamma^*$, $1\leq\alpha<\infty$ and $1\leq\beta\leq\infty$, the solution $(\varphi,
  \eta)$ of the system $(2.19)$ satisfies the following estimates:
\begin{equation}
  \|\varphi(\cdot,t)-\varphi_{\infty}\|_{X_{\alpha\beta}}
  +\|\eta(t)-\eta_{\infty}\|_{Y_{\beta}}\leq
  C_{\alpha\beta}e^{-\lambda t}[\|\varphi_0\|_{X_{\alpha\beta}}
  +\|\eta_0\|_{Y_{\beta}}] \quad \mbox{for}\;\; t\geq 0,
\end{equation}
  where $\varphi_0=\varphi(\cdot,0)$ and $\eta_0=\eta(0)$ are the initial data of
  $\varphi$ and $\eta$, respectively, $\varphi_{\infty}$ and $\eta_{\infty}$ are
  functions in $\mathbb{B}(0,R_s)$ and $\mathbb{S}^{n-1}$, respectively, having
  the following expressions:
\begin{equation}
  \varphi_{\infty}(x)=-p_s'(r)\sum_{l=1}^nc_lY_{1l}(\omega), \qquad
  \eta_{\infty}(\omega)=\sum_{l=1}^nc_lY_{1l}(\omega),
\end{equation}
  where $c_1,c_2,\cdots,c_n$ are real constants uniquely determined by the
  initial data $\varphi_0$ and $\eta_0$, and $C_{\alpha}$ is a positive constant.
  Moreover, in case $\alpha=\infty$ we have also the following estimate:}
\begin{equation}
  \|\varphi(\cdot,t)-\varphi_{\infty}\|_{\dot{X}_{\infty\beta}}
  +\|\eta(t)-\eta_{\infty}\|_{Y_{\beta}}\leq
  C_{\beta}e^{-\lambda t}[\|\varphi_0\|_{\dot{X}_{\infty\beta}}
  +\|\eta_0\|_{Y_{\beta}}] \quad \mbox{for}\;\; t\geq 0.
\end{equation}

  {\em Proof}:\ \ This is an immediate consequence of Lemmas 7.1--7.3.  $\quad\Box$
\medskip

  Hence, we have finished proving Theorem 1.1.  $\quad\Box$


\begin{thebibliography}{99}
\bibitem{tumrev1} R. P. Araujo and D. L. McElwain, A history of the study of
  solid tumor growth: the contribution of mathematical modeling, \textit{Bull.
  Math. Biol.}, {66}(2004), 1039--1091.
\bibitem{ChenFri}  X. Chen and A. Friedman,  A free boundary problem for an
  elliptic-hyperbolic system: an application to tumor growth, \textit{SIAM J.
  Math. Anal.}, {35}(2003), 974--986.
\bibitem{ChenCuiF} X. Chen, S. Cui and A. Friedman, A hyperbolic free boundary
  problem modeling tumor growth: asymptotic behavior. \textit{Trans. Amer.
  Math. Soc.}, {357}(2005), 4771-4804.
\bibitem{Cui1} S. Cui, Existence of a stationary solution for the modified
  Ward-King tumor growth model, \textit{Advances in Appl. Math.},
  {36}(2006), 421--445.
\bibitem{Cui2} S. Cui, Asymptotic stability of the stationary solution for
  a hyperbolic free boundary problem modeling tumor growth, \textit{SIAM
  J. Math. Anal.}, {40}(2008), 1692--1724.
\bibitem{Cui3} S. Cui,  Lie group action and stability analysis of
  stationary solutions for a free boundary problem modelling tumor growth,
  \textit{J. Diff. Equa.}, {246}(2009), 1845--1882.
\bibitem{Cui4} S. Cui, Asymptotic stability of the stationary solution for
  a parabolic-hyperbolic free boundary problem modeling tumor growth,
  \textit{SIAM J. Math. Anal.}, {45}(2013), proofs.
\bibitem{CuiEsc1} S. Cui and J. Escher, Bifurcation analysis of an elliptic
  free boundary problem modelling the growth of avascular tumors,
  \textit{SIAM J. Math. Anal.}, {39}(2007), 210--235.
\bibitem{CuiEsc2} S. Cui and J. Escher, Asymptotic Behavior of Solutions of
  Multidimensional Moving Boundary Problem Modeling Tumor Growth,
  {\em Comm. Part. Diff. Equa.}, {33}(2008), 636--655.
\bibitem{CuiEsc3} S. Cui and J. Escher, Well-posedness and stability of
  a multidimentional moving boundary problem modeling the growth of tumors,
  \textit{Arch. Rat. Mech. Anal.}, {191}(2009), 173--193.
\bibitem{CuiFri1} S. Cui and A. Friedman, A free boundary problem for a
  singular system of differential equations: an application to a model of
  tumor growth, \textit{Trans. Amer. Math. Soc.}, {355}(2002),
  3537--3590.
\bibitem{CuiFri2} S. Cui and A. Friedman, A hyperbolic free boundary problem
  modeling tumor growth, \textit{Interfaces and Free Bound.},
  {5}(2003), 159-181.
\bibitem{CuiWei} S. Cui and X. Wei, Global existence for a parabolic-hyperbolic
  free boundary problem modeling tumor growth, \textit{Acta Math. Appl. Sinica
  English Series}, {21}(2005), 597--614.
\bibitem{EN}  K. J. Engel and R. Nagel, {\em One-Parameter Semigroups for
  Linear Evolution Equations}, Springer-verlag, New York, 2000.
\bibitem{Esc} J. Escher, Classical solutions to a moving boundary problem
  for an elliptic-parabolic system, \textit{Interfaces Free and Boundaries},
  {6} (2004), 175--193.
\bibitem{EscSim} J. Escher and G. Simonett, Classical solutions for Hele-Shaw
  models with surface tension, \textit{Adv. Diff. Equa.}, {2} (1997),
  619--642.
\bibitem{tumrev2} A. Fasano, A. Bertuzzi and A. Gandolfi, Mathematical
  modelling of tumour growth and treatment, \textit{Lect. Notes Math.},
  {1872}(2006), 71--106.
\bibitem{FonFri1} M. A. Fontelos and A. Friedman, Symmetry-breaking bifurcations
  of free boundary problems in three dimensions, \textit{Asymp. Anal.},
  {35}(2003), 187--206.
\bibitem{Fried1} A. Friedman,  A hierarchy of cancer models and their
  mathematical challenges, \textit{Disc. Cont. Dyna. Syst.} (B),
  {4}(2004),  147--159.
\bibitem{Fried2} A. Friedman: Cancer models and their mathematical analysis,
  \textit{Lect. Notes Math.}, {1872}(2006), 223--246.
\bibitem{Fried3} A. Friedman,  Mathematical analysis and challenges arising
  from models of tumor growth, \textit{Math. Model. Meth. Appl. Sci.},
  {17}(2007), suppl., 1751--1772.
\bibitem{FriHu1} A. Friedman and B. Hu, Bifurcation from stability to
  instability for a free boundary problem arising in tumor model,
  \textit{Archive Rat. Mech. Anal.}, {180}(2006), 293--330.
\bibitem{FriHu2} A. Friedman and B. Hu, Asymptotic stability for a free
  boundary problem arising in a tumor model, \textit{J. Diff. Egua.},
  {227}(2006), 598--639.
\bibitem{FriHu3} A. Friedman and B. Hu, Stability and instability of
  Liapounov-Schmidt and Hopf bifurcation for a free boundary problem arising
  in a tumor model, \textit{Trans. Amer. Math. Soc.},
\bibitem{FriRei1} A. Friedman and F. Reitich, Analysis of a mathematical
  model for growth of tumors, \textit{J. Math. Biol.}, {38}(1999),
  262--284.
\bibitem{FriRei2} A. Friedman and F. Reitich, Symmetry-breaking bifurcation
  of analytic solutions to free boundary problems: An application to a model
  of tumor growth, \textit{Trans. Amer. Math. Soc.}, {353}(2000),
  1587--1634.
\bibitem{Green1} H. P. Greenspan, Models for the growth of solid tumor by diffusion,
  \textit{Stud. Appl. Math.} {51}(1972), 317--340.
\bibitem{Green2} H. P. Greenspan, On the growth and stability of cell cultures
  and solid tumors, \textit{J. Theor. Biol.}, {56}(1976), 229--242.
\bibitem{tumrev3} J. S. Lowengrub, H. B. Frieboes, F. Jin et al, Nonlinear
  modelling of cancer: bridging the gap between cells and tumours,
  \textit{Nonlinearity}, {23}(2010), R1--R91.
\bibitem{Pazy} A. Pazy, \textit{Semigroups of Linear Operators and
  Applications to Partial Differential Equations}, Springer, New York: 1983.
\bibitem{PPM} G. Pettet, C. Please and M. McElwain, The migration of
  cells in multicell tumor  spheroids, \textit{Bull. Math. Biol.},
  {63}(2001), 231--257.
\bibitem{TP} M. J. Tindall and C. P. Please, Modelling the cell cycle and
  cell movement in multicellular tumour spheroids, \textit{Bull. Math. Biol.},
  {69}(2007), 1147--1165.
\bibitem{SW} E. M. Stein and G. Weiss, \textit{Introduction to Fourier
  Analysis on Euclidean Spaces, Chapter IV}, Princeton University Press, New
  Jersey: 1971.
\bibitem{WuCui} J. Wu and S. Cui, Asymptotic stability of stationary solutions
  of a free boundary problem modelling the growth of tumors with fluid tissues,
  \textit{SIAM J. Math. Anal.}, {41}(2009), 391--414.
\bibitem{WZ} J. Wu and F. Zhou, Asymptotic behavior of solutions of a free
  boundary problem modelling the growth of tumors with fluid-like tissue under
  the action of inhibitors, \textit{Trans. Amer. Math. Soc.}, {365}(2013),
  4181--4207.
\end{thebibliography}
\end{document}